\documentclass{amsart}
\usepackage{amsmath,amsxtra,amssymb}
\usepackage[all]{xy}

\newcommand{\hz}{\vspace{0.5cm}}

\renewcommand{\qed}{\hspace*{\fill}$\Box$\hz\pagebreak[1]}

\newcommand{\B}{{\mathcal B}}

\newcommand{\M}{{\mathcal M}}

\newcommand{\N}{{\mathcal N}}

\newtheorem{theorem}{Theorem}[section]
\newtheorem{cor}[theorem]{Corollary}
\newtheorem{prop}[theorem]{Proposition}
\newtheorem{rem}[theorem]{Remark}
\newtheorem{lemma}[theorem]{Lemma}
\newtheorem{defi}[theorem]{Definition}
\newtheorem{exam}[theorem]{Example}

\begin{document}
\title{Classification of hyperfinite factors up to completely bounded isomorphism of their preduals}
\author{Uffe Haagerup$^{(1)}$ and Magdalena Musat$^{(2)}$}
\address{$^{(1)}$ Department of Mathematics and
Computer Science, University of Southern Denmark, Campusvej 55, 5230
Odense M, Denmark.\\
$^{(2)}$Department of Mathematics and Computer Science, University
of Southern Denmark, Campusvej 55, 5230 Odense M, Denmark and
Department of Mathematical Sciences, University of Memphis, 373 Dunn
Hall, Memphis, TN, 38152, USA.}
\email{$^{(1)}$haagerup@imada.sdu.dk\\$^{(2)}$musat@imada.sdu.dk,
mmusat@memphis.edu}.

\date{}


\maketitle

\begin{center}
{\em Dedicated to Alain Connes on his 60th birthday}
\end{center}

\begin{abstract}
In this paper we consider the following problem: When are the
preduals of two hyperfinite (=injective) factors $\M$ and $\N$ (on
separable Hilbert spaces) cb-isomorphic (i.e., isomorphic as
operator spaces)? We show that if $\M$ is semifinite and $\N$ is
type III, then their preduals are not cb-isomorphic. Moreover, we
construct a one-parameter family of hyperfinite type III$_0$-factors
with mutually non cb-isomorphic preduals, and we give a
characterization of those hyperfinite factors $\M$ whose preduals
are cb-isomorphic to the predual of the unique hyperfinite type
III$_1$-factor. In contrast, Christensen and Sinclair proved in 1989
that all infinite dimensional hyperfinite factors with separable
preduals are cb-isomorphic and more recently, Rosenthal, Sukochev
and the first-named author proved that all hyperfinite type
III$_\lambda$-factors, where $0< \lambda\leq 1$, have cb-isomorphic
preduals.
\end{abstract}

\section{Introduction and formulation of the main results}
\setcounter{equation}{0}

In the paper \cite{ChS1}, Christensen and Sinclair proved that if
$\M$ and $\N$ are infinite dimensional factors with separable
preduals, then $\M$ and $\N$ are cb-isomorphic
($\M\overset{\text{cb}}{\simeq} \N$)\,, i.e., there exists a linear
bijection $\phi$ of $\M$ onto $\N$ such that both $\phi$ and
$\phi^{-1}$ are completely bounded. In 1993 Kirchberg (cf.
\cite{Ki}) proved a similar result for $C^*$-algebras, namely if $A$
and $B$ are simple, separable, nuclear, non-type I $C^*$-algebras,
then $A\overset{\text{cb}}{\simeq} B$\,.

However, if one turns to preduals of von Neumann algebras (on
separable Hilbert spaces), the situation is very different.
Rosenthal, Sukochev and the first-named author proved in \cite{HRS}
that if $\M$ is a II$_1$-factor and $\N$ is a properly infinite von
Neumann algebra, then their preduals $\M_*$ and $\N_*$ are not
isomorphic as Banach spaces, so in particular they are not
cb-isomorphic. Moreover, the Banach space isomorphism classes
(respectively, cb-isomorphism classes) of separable preduals of
hyperfinite and semifinite von Neumann algebras are completely
determined by \cite{HRS}, Theorem 5.1.

By a combination of two recent results of Pisier and Junge, the
predual $\M_*$ of a semifinite factor $\M$ cannot be cb-isomorphic
to the predual of $R_\infty$\,, the unique hyperfinite factor of
type III$_1$\,, because Pisier's operator Hilbert space $OH$ does
not cb-embed in $\M_*$ by \cite{Pi3}, while $OH$ admits a
cb-embedding into $(R_\infty)_*$, as proved in \cite{Ju} (see also
\cite{HM})\,.

A von Neumann algebra $\M$ with separable predual is hyperfinite (or
approximative finite dimensional) if $\M$ is the strong closure of
an increasing union $\cup_{n=1}^\infty \M_n$ of finite dimensional
$*$-subalgebras $\M_n$\,. By Connes' celebrated work \cite{Con}\,, a
factor $\M$ (with separable predual) is hyperfinite if and only if
it is injective. It is well-known that the same holds for
non-factors (see Section 6 of \cite{Ha5} and the references given
therein), so in the following we do not need to distinguish between
"hyperfiniteness" and "injectivity" for von Neumann algebras with
separable preduals.

The main results of this paper are the following:
\begin{theorem}\label{th1}
Let $\M$ and $\N$ be hyperfinite von Neumann algebras with separable
preduals $\M_*$ and $\N_*$, respectively. If $\M$ is type III and
$\N$ is semifinite, then the preduals $\M_*$ and $\N_*$ are not
cb-isomorphic. More generally, $\M_*$ is not cb-isomorphic to a
cb-complemented subspace of $\N_*$\,.
\end{theorem}

\begin{theorem}\label{th2}
The predual of a hyperfinite type III-factor $\M$ (on a separable
Hilbert space) is cb-isomorphic to the predual of the hyperfinite
type III$_1$-factor $R_\infty$ if and only if there exists a normal
invariant state on the flow of weights $(Z(\N), (\theta_s)_{s\in
\mathbb{R}})$ for $\M$\,.
\end{theorem}

\begin{theorem}\label{th3}
There exist uncountably many cb-isomorphism classes of preduals of
hyperfinite type III$_0$-factors (on separable Hilbert spaces).
\end{theorem}

Both Theorems \ref{th2} and \ref{th3} rely on Connes' classification
of type III-factors (cf. \cite{Co4}), and the Connes-Takesaki flow
of weights for type III-factors (cf. \cite{CoT}) and on Connes'
classification of injective factors (see \cite{Con}), which we will
outline below.

In \cite{Co4} Connes introduced  the subclassification of type III-
factors into type III$_\lambda$-factors, where $0\leq \lambda\leq
1$\,. Later, in \cite{CoT} Connes and Takesaki introduced the
"smooth flow of weights" (now called "the flow of weights") of a
type III-factor. Following Takesaki's exposition in \cite{Tk}, Vol.
II, pp. 364-368, the flow of weights can be constructed as follows.
Let $\M$ be a type III-factor, and let $\phi$ be a normal, faithful
state on $\M$\,. Consider the crossed product
$\N:=\M\rtimes_{\sigma^\phi} \mathbb{R}$\,, where
$(\sigma_t^\phi)_{t\in \mathbb{R}}$ is the modular automorphism
group associated with $\phi$\,. Then $\N$ is generated by an
embedding $\pi(\M)$ of $\M$ into $\N$ and by a one-parameter group
$(\lambda(t))_{t\in \mathbb{R}}$ of unitaries in $\N$\,. Moreover,
there is a s.o.t-continuous dual action
$(\widetilde{\theta_s})_{s\in \mathbb{R}}$ of $\mathbb{R}$ on
$\N$\,, characterized by the relations
\begin{eqnarray*}
\widetilde{\theta}_s(\pi(x))&=&\pi(x)\,, \quad x\in \M\\
\widetilde{\theta}_s(\lambda(t))&=&e^{ist} \lambda(t)\,, \quad t\in
\mathbb{R}\,,
\end{eqnarray*}
for all $s\in \mathbb{R}$\,. Let $\theta_s$ be the restriction of
$\widetilde{\theta}_s$ to the center $Z(\N)$ of $\N$\,. Then
$(Z(\N), (\theta_s)_{s\in \mathbb{R}})$ is called {\em the flow of
weights for} $\M$\,. It is independent (up to isomorphism) of the
choice of the state $\phi$ on $\M$\,. Since $\M$ is a factor,
$\theta=({\theta}_s)_{s\in \mathbb{R}}$ acts ergodically on
$Z(\N)$\,, i.e., the fixed point algebra $Z(\N)^\theta$ for the
action $\theta$ is equal to $\mathbb{C}1$\,. Since $Z(\N)\simeq
L^\infty(\Omega, \mu)$ for some standard Borel measure space
$(\Omega, \mu)$\,, the flow $\theta=(\theta_s)_{s\in \mathbb{R}}$
can be realized as the flow associated to a one-parameter family
$(\sigma_s)_{s\in \mathbb{R}}$ on non-singular Borel transformations
of $(\Omega, \mu)$\,, that is, for all $s\in \mathbb{R}$\,,
\begin{eqnarray*}
(\theta_s f)(x)&=& f(\sigma_s^{-1}(x))\,, \quad f\in
L^\infty(\Omega, \mu)\,, x\in \Omega\,.
\end{eqnarray*}
The connection between Connes' type III$_\lambda$-factors and the
flow of weights is given by {\bf 1.4.}, {\bf 1.5.} and {\bf 1.6.}
below (cf. \cite{CoT}\,, \cite{Ta2} and \cite{Tk}, Vol. II, Chapter
XII).

Let $\M$ be a type III-factor (with separable predual).
Then\\[0.2cm]
\noindent {\bf 1.4.} $\M$ is of type III$_0$ if and only if the flow
$(\theta_s)_{s\in \mathbb{R}}$ is non-periodic, i.e., $\theta_s\neq
\text{Id}_{Z(\N)}$, for all $s\in \mathbb{R}\smallsetminus\{0\}$\,.
In this case the flow is non-transitive (=properly ergodic), which
means that the measure $\mu$ described above is not concentrated on
a single $\sigma$-orbit in $\Omega$\,.\\[0.3cm]
\noindent {\bf 1.5.} For $0< \lambda< 1$\,, $\M$ is of type
III$_\lambda$ if and only if the flow $(\theta_s)_{s\in \mathbb{R}}$
is periodic with minimal period equal to $-\log{\lambda}$\,. In this
case, \[ Z(\N)\simeq
L^\infty(\mathbb{R}/{((-\log{\lambda})\mathbb{Z})})\,.
\] and for all $s\in \mathbb{R}$\,, $\theta_s$ is induced by the
translation $\sigma_s: x\mapsto x+s$ on
$\Omega=\mathbb{R}/{((-\log{\lambda})\mathbb{Z})}$\,.\\[0.3cm]
\noindent {\bf 1.6.} $\M$ is of type III$_1$ if and only if
$\theta_s=\text{Id}_{Z(\N)}$\,, for all $s\in \mathbb{R}$\,. In this
case $Z(\N)=\mathbb{C}1$\,.

\vspace*{0.3cm} For hyperfinite type III-factors with separable
predual much more is known, owing to Connes' classification of
injective factors \cite{Con}, and related work by Krieger
\cite{Kr1}, Connes \cite{Co7} and by the
first-named author \cite{Ha2}, namely:\\[0.2cm]
\noindent {\bf 1.7.} The map $\M\mapsto ((Z(\N), (\theta_s)_{s\in
\mathbb{R}})$ gives a one-to one correspondence between the set of
(isomorphism classes of) hyperfinite type III$_0$-factors onto the
set of (isomorphism classes of) non-transitive ergodic flows $(A,
(\theta_s)_{s\in \mathbb{R}})$ on abelian von Neumann algebras $A$
with separable predual. In particular, there are uncountably many
isomorphism classes of hyperfinite type III$_0$-factors (cf.
\cite{Con}, \cite{Kr1}).\\[0.3cm]
\noindent {\bf 1.8.} For each $0< \lambda< 1$, there is exactly one
(up to isomorphism)  hyperfinite factor of type III$_\lambda$\,,
namely the Powers factor \begin{equation*}
R_\lambda:=\otimes_{n=1}^\infty (M_2(\mathbb{C}), \phi_\lambda)\,,
\end{equation*} where
$\phi_\lambda=\text{Tr}(h_\lambda\,\cdot\,)$\,, $\text{Tr}$ being
the non-normalized trace on $M_2(\mathbb{C})$ and
$h_\lambda={\frac1{1+\lambda}}\left(
\begin{array}
[c]{cc}%
\lambda & 0\\
0 & 1
\end{array}
\right)$\,. (See \cite{Con}.)\\[0.2cm]
\noindent {\bf 1.9.} There is only one (up to isomorphism)
hyperfinite type III$_1$-factor, namely the Araki-Woods factor
$R_\infty$\,, which can be expressed as the (von Neumann algebra)
tensor product $R_\infty=R_{\lambda_1}\bar{\otimes} R_{\lambda_2}$
of two Powers factors with
${\frac{\log{\lambda_1}}{\log{\lambda_2}}}\notin \mathbb{Q}$ (cf.
\cite{Co7} and \cite{Ha2}).\\[0.2cm]
Proofs of ${\bf 1.7.}$\,, ${\bf 1.8.}$ and ${\bf 1.9.}$ can also be
found in \cite{Tk}, Vol. III, Chap. XVIII.

The rest of the paper is organized in the following way. In Section
2, we obtain some Stinespring-type decomposition results for
completely positive and completely  bounded maps, which will allow
us to show that two properly infinite hyperfinite von Neumann
algebras $\M$ and $\N$ with  separable preduals have cb-isomorphic
preduals if and only if there exist von Neumann algebras embeddings
$i: \M\rightarrow \N$\,, $j:\N\rightarrow \M$ and normal conditional
expectations $E:\N{\rightarrow} i(\M)$\,, $F:\M{\rightarrow}
j(\N)$\,. From this result, Theorem \ref{th1} follows easily by
results of Sakai \cite{Sak2} and Tomiyama \cite{To} on normal
conditional expectations.

In Section 3 we prove Theorem \ref{th2}. The most difficult part is
to show that $R_\infty$ embeds into $\M$ as the range of a normal
conditional expectation, provided that there exists a normal,
invariant state on the flow of weights for $\M$\,. This part of the
proof relies heavily on the main result from \cite{HS} on the
classification of normal states on a von Neumann algebra up to
approximative unitary equivalence.

Finally, in Section 4 we prove Theorem \ref{th3} by giving an
explicit construction of a one-parameter family $(A^{(t)}\,,
\theta^{(t)})_{0\leq t< 2}$ of non-transitive, ergodic flows
$\theta^{(t)}=(\theta_s^{(t)})_{s\in \mathbb{R}}$ on abelian von
Neumann algebras $A^{(t)}$ with separable preduals $A_*^{(t)}$\,,
satisfying the following property:
\begin{equation}\label{eq2211221144224433}
\lim\limits_{n\rightarrow \infty}\|\omega\circ
\theta_{2^n}^{(t)}-\omega\|=t\,, \quad \omega\in A_*^{(t)}\,.
\end{equation}
Then by ${\bf 1.7.}$ above, $(A^{(t)}\,, \theta^{(t)})_{0\leq t< 2}$
are the flows of weights associated with hyperfinite type III$_0$-
factors $(\M^{(t)})_{0\leq t< 2}$\,, and by
(\ref{eq2211221144224433}) combined with the results of Section 2,
we obtain that ${\M_*^{(t_1)}}$ and ${\M_*^{(t_2)}}$ are not
cb-isomorphic when $t_1\neq t_2$\,. It is interesting to note that
the factors $(\M^{(t)})_{0\leq t< 2}$ cannot be separated by Connes'
$S$- and $T$-invariants. Being type III$_0$-factors,
$S(\M^{(t)})=\{0, 1\}$\,, for all $t\in [0, 2)$, and in Theorem
\ref{th9322} we prove that $T(\M^{(t)})=\left\{\frac{2\pi k}{2^n};
k\in \mathbb{Z}\,, n\in \mathbb{N}\right\}$\,, for all $t\in [0,
2)$\,.

For details on operator spaces and completely bounded maps we refer
to the monographs \cite{ER, Pi1}. We shall briefly recall some
definitions that are relevant for our paper. An operator space $X$
is a Banach space given together with an isometric embedding
$X\subseteq {\mathcal B}(H)$\,, the algebra of bounded linear
operators on a Hilbert space $H$\,. For all $n\geq 1$, this
embedding determines a norm on $M_n(X)$ (the $n\times n$ matrices
over $X$), induced by the space $M_n({\mathcal B}(H))\cong {\mathcal
B}(H^n)$\,. The morphisms in the category of operator spaces are
{\em completely bounded maps}. Given a linear map $\phi :
X\rightarrow Y$ between two operator spaces $X$ and $Y$ and $n\geq
1$\,, define ${\phi}_n : {M_n(X)} \rightarrow {M_n(Y)}$ by
${\phi}_n([x_{ij}])= [\phi(x_{ij})]$\,, for all $
[x_{ij}]_{i,j=1}^n\in M_n(X)$\,. Let $\|\phi\|_{cb}:=\sup
\{\|{\phi}_n\|\,; n\in {\mathbb{N}} \,\}$\,. The map $\phi$ is
called {\em completely bounded} (for short, {\em cb}) if
$\|\phi\|_{cb} <\infty\,,$ and $\phi$ is called {\em completely
isometric} if all ${\phi}_n$ are isometries. The space of all cb
maps from $X$ to $Y$\,, denoted by $\mathcal{C}\B(X, Y)$\,, is an
operator space with matrix norms defined by $M_n(\mathcal{C}\B(X,
Y))=\mathcal{C}\B(X, M_n(Y))$\,, $n\geq 1$\,. The dual of an
operator space $X$ is, again, an operator space
$X^*=\mathcal{C}\B(X, \mathbb{C})$\,. A von Neumann algebra $\M$
carries a natural operator space structure, and its predual $\M_*$
carries the operator space structure induced by the completely
isometric embedding into the dual $\M^*$ of $\M$\,.

\section{Stinespring-type decomposition theorems and applications}
\setcounter{equation}{0}

\begin{lemma}\label{lem6789}
Let $\M$ and $\N$ be von Neumann algebras with separable preduals
$\M_*$ and $\N_*$\,, respectively. Let $\alpha: \M\rightarrow \N$ be
a completely positive map. Then there exists a completely positive
unital map $\widetilde{\alpha}: \M\rightarrow \N$ such that
\[ \alpha(a)=\alpha(1)^{1/2}\widetilde{\alpha}(a)\alpha(1)^{1/2}\,, \quad a\in
\M\,,
\] where $1$ denotes the unit of $\M$\,. Moreover, if $\alpha$ is normal,
then $\widetilde{\alpha}$ can be chosen to be normal.
\end{lemma}

\begin{proof}
Assume first that $\text{supp}(\alpha(1))=1_{\N}$\,, where $1_{\N}$
denotes the unit of $\N$\,. Let $H$ be a separable Hilbert space
with $\N\subseteq {\mathcal B}(H)$\,. The operator $\alpha(1)^{1/2}$
is one-to-one and has dense range, denoted by $H_0$\,. Hence we get
a (possibly unbounded) map $\alpha(1)^{-1/2}: H_0\rightarrow H$\,.

Now fix $a\in \M$\,, $a\geq 0$\,, and define a positive sesquilinear
form on $H_0\times H_0$ by
\[ s(x, y):=\langle \alpha(a)\alpha(1)^{-1/2}x,
\alpha(1)^{-1/2}y\rangle\,, \quad x, y\in H_0\,. \] Note that $s$ is
positive, since $\alpha$ is so. We now show that $s$ is a bounded
sesquilinear form. For all $x\in H_0$\,,
\begin{equation*}
s(x, x)\leq \|a\|\langle \alpha(1)\alpha(1)^{-1/2}x,
\alpha(1)^{-1/2}x\rangle=\|a\|\|x\|^2\,.
\end{equation*}
By Schwarz's inequality,
\[ |s(x, y)|\leq s(x, x)^{1/2}s(y, y)^{1/2}\leq \|a\|\|x\|\|y\|\,,
\quad x, y\in H_0\,. \] Hence there exists a unique operator $T\in
{\mathcal B}(H)$ such that
\[ \langle Tx, y\rangle=\langle \alpha(a)\alpha(1)^{-1/2}x\,,
\alpha( 1)^{-1/2}y\rangle\,, \quad x, y\in H_0\,. \] Note first that
$T\in \N$\,. This follows from the fact that for all $x, y\in H_0$
and all unitaries $U$ in the commutant ${\N}'$ of $\N$\,,
\[ \langle TUx\,, Uy\rangle=\langle Tx\,, y\rangle\,, \] wherein we
use the fact that $\alpha(a)\in \N$, and $\alpha(1)^{-1/2}$ is
affiliated with $\N$\,. Clearly $T$ is positive.

Since $\M$ is the span of its positive part $\M_{+}$\,, we infer
that for all $a\in \M$\,, there exists a unique element
$\widetilde{\alpha}(a)\in \N$ such that
\[ \langle \widetilde{\alpha}(a)x, y\rangle=\langle
\alpha(a)\alpha(1)^{-1/2}x\,, \alpha(1)^{-1/2}y\rangle\,, \quad x,
y\in H_0\,. \] By uniqueness, the map $\widetilde{\alpha}$ is
linear. Also, clearly $\widetilde{\alpha}(a)\geq 0$\,, whenever
$a\geq 0$\,, and $\widetilde{\alpha}(1)=1_{\N}$\,. Looking at
$n\times n$ matrices over $\M$ we infer that $\widetilde{\alpha}$ is
completely positive. Moreover, for all $a\in \M$\,,
\[
\alpha(a)=\alpha(1)^{1/2}\widetilde{\alpha}(a)\alpha(1)^{1/2}\,,
\] since $\langle \alpha(a)x, y\rangle = \langle
\widetilde{\alpha}(a)\alpha(1)^{1/2}x, \alpha(1)^{1/2}y\rangle$\,,
for all $x, y\in H$\,. Note also that $\widetilde{\alpha}$ is normal
if $\alpha$ is so.

It remains to consider the case when
$p:=\text{supp}(\alpha(1))\lvertneqq 1_{\N}$\,. Apply the previous
argument to the mapping $\alpha: \M\rightarrow p{\N}p\subseteq
\N$\,. We then obtain a completely positive map $\widetilde{\alpha}:
\M\rightarrow p{\N}p$ such that $\widetilde{\alpha}(1)=p$ and
$\alpha(a)=\alpha(1)^{1/2}\widetilde{\alpha}(a)\alpha(1)^{1/2}$\,,
for all $a\in \M$\,. Choose a normal state $\phi$ on $\M$\,, and set
\[ \widehat{\alpha}(a):=\widetilde{\alpha}(a)+\phi(a)(1_{\N}-p)\,,
\quad a\in \M\,. \] Then $\widehat{\alpha}$ is completely positive
and $\widehat{\alpha}(1)=p+(1_{\N}-p)=1_{\N}$\,. Furthermore, since
$\text{supp}({\alpha(1)}^{1/2})=p$\,,
\[
\alpha(1)^{1/2}\widehat{\alpha}(a)\alpha(1)^{1/2}=\alpha(a)+\phi(a)\alpha(1)^{1/2}(1_{\N}-p)\alpha(1)^{1/2}=\alpha(a)\,.
\]
Moreover, if $\alpha$ is normal, then both $\widetilde{\alpha}$ and
$\widehat{\alpha}$ are normal. The proof is complete.
\end{proof}

\begin{lemma}\label{lem82988298}
Let $\N$ be a von Neumann algebra with separable predual. If
$p\in\N$ is a properly infinite projection with central support
equal to the identity $1$ of $\N$\,, then $p\sim 1$\,.
\end{lemma}
\begin{proof}
This result is well-known and it follows by standard comparison
theory of projections. For convenience of the reader, we include a
proof. Assume that $p\neq 1$. Choose a maximal family $(p_i)_{i\in
I}$ of pairwise orthogonal non-zero projections such that $p_i\prec
p$\,, for all $i\in I$\,. Note that $I$ must be countable. We first
show that
\begin{equation}\label{eq4444333322228}\sum\limits_{i\in I}
p_i=1\,. \end{equation} Suppose by contradiction that
$\sum\limits_{i\in I} p_i< 1$\,. Set $q:=1-\sum\limits_{i\in I}
p_i$\,. If $c(p)$ and $c(q)$ denote the central support of $p$ and
$q$, respectively, then $c(q)c(p)=c(q)\neq 0$\,. This implies that
there exist nonzero projections $q_0$ and $p_0$ such that $q_0\leq
q$\,, $p_0\leq p$ and $q_0\sim p_0$\,. Hence $q_0\prec p$\,.  Since
$q_0\leq q$\,, this contradicts the maximality assumption of the
family $(p_i)_{i\in I}$\,, and (\ref{eq4444333322228}) is proved.
Since $p$ is properly infinite, we can write
\begin{equation}\label{eq4444333322229} p=\sum\limits_{i\in
I} r_i\,, \end{equation} where $(r_i)_{i\in I}$ are pairwise
orthogonal projections so that $p\sim r_i$\,, for all $i\in I$.
Hence $p_i\prec p\sim r_i$\,, for all $i\in I$\,. Together with
(\ref{eq4444333322228}) and (\ref{eq4444333322229}), this implies
that $1\preceq p$\,. Clearly $p\preceq 1$\,, and therefore $p\sim
1$\,.
\end{proof}
We now prove the following Stinespring-Kasparov-type theorem (see
\cite{Kas}, Theorem 3(1)):

\begin{theorem}\label{th33}
Let $\M$ and $\N$ be von Neumann algebras with separable preduals
$\M_*$ and $\N_*$\,, respectively. Assume, moreover, that $\N$ is
properly infinite. Let $\alpha: \M\rightarrow \N$ be a normal,
completely positive map. Then there exists a normal
$*$-representation $\pi: \M\rightarrow \N$ and an operator $V\in \N$
such that
\[ \alpha(a)=V^*\pi(a)V\,, \quad a\in \M\,. \]
\end{theorem}

\begin{proof}
By Lemma \ref{lem6789} we can assume without loss of generality that
$\alpha(1)=1_{\N}$\,, where $1$ and $1_{\N}$ are the identities of
$\M$ and $\N$, respectively. Following Stinespring's construction
(see \cite{Pau}, Theorem 4.1), we define a positive sesquilinear
form $s$ on the algebraic tensor product $\M\odot H$\,, where $H$ is
a separable Hilbert space with $\N\subseteq {\mathcal B}(H)$\,, by
\[ s(a\otimes x, b\otimes y):=\langle \alpha(b^*a)x, y\rangle\,,
\quad a, b\in \M\,, x, y\in H\,. \] Let $L:=\{z\in \M\odot H; s(z,
z)=0\}$\,, and note that $({\M\odot H})/L$ is a prehilbert space
whose completion we denote by $K$\,. For all $a\in \M$ and $x\in H$
let $[a\otimes x]$ denote the corresponding element in the quotient
space ${{\M\odot H}/L}$\,. For all $a\in \M$ define $\pi_0(a)$ by
\[ \pi_0(a)[b\otimes x]:=[ab\otimes x]\,, \quad b\in \M\,,
x\in H\,. \] Then $\pi_0(a)$ is a densely-defined, bounded operator
on the dense subspace $({\M\odot H})/L$ of $K$\,. Hence $\pi_0(a)$
extends to a bounded linear operator on the whole $K$\,, and this
yields a map $\pi_0: \M\rightarrow {\mathcal B}(K)$\,. It is easily
checked that $\pi_0$ is a unital $*$-representation. Moreover,
$\pi_0$ is normal. This follows immediately from the fact that
$\alpha$ is normal and that
\begin{equation}\label{eq9999999945327}
\langle \pi_0(a)[b\otimes x]\,, [c\otimes y]\rangle=\langle
\alpha(c^*ab)x\,, y\rangle\,, \end{equation} for all $a, b, c\in
\M$\,, and all $x, y\in H$\,. Define now $W: H\rightarrow K$ by\[
Wx:=[1\otimes x]\,, \quad x\in H\,. \] By (\ref{eq9999999945327}),
it follows that $\langle \pi_0(a)Wx, Wy\rangle=\langle \alpha(a)x,
y\rangle$\,, for all $a\in \M$ and all $x, y\in H$\,. Since
$\alpha(1)=1_{\N}$\,, $W$ is an isometry of $H$ into $K$\,,
$W^*W=1_{\N}$ and
\begin{equation}\label{eq5555555344445}
\alpha(a)=W^*\pi_0(a)W\,, \quad a\in \M\,.
\end{equation}
Let $\N^{\prime}$ denote the commutant of $\N$ in ${\mathcal
B}(H)$\,. Define a normal $*$-representation of $\N^{\prime}$ on
$({\M\odot H})/L$ by
\[ \sigma_0(n')[a\otimes x]:=[a\otimes n'x]\,, \quad n'\in \N^{\prime}\,. \]
It is easily checked that $\sigma_0$ is a well-defined linear map on
$({\M\odot H})/L$\,. Now let $n'\in {\N'}$\,. Then
\begin{eqnarray*}
\|[a\otimes n'x]\|^2=\langle \alpha(a^*a)n'x, n'x\rangle =\langle
(n')^*\alpha(a^*a)n'x, x\rangle &=&\langle (n')^*n'\alpha(a^*a)x,
x\rangle\\&\leq &\|n'\|^2\|[a\otimes x]\|^2\,,
\end{eqnarray*}
wherein we have used the fact that $\alpha(\M)\subseteq \N$\,. We
deduce that $\sigma_0(n')$ has a unique extension to an operator
$\sigma(n')\in {\mathcal B}(K)$ such that
\begin{equation}\label{eq7777644443233} \sigma(n')[a\otimes x]=[a\otimes n'x]\,,
\quad a\in \M\,, x\in H\,.
\end{equation} Note also that $\sigma(1_{\N})=1_K$\,, and it is easily checked
that the map $\sigma:\N'\rightarrow {\mathcal B}(K)$ thus defined is
a normal $*$-representation. Set
\[ \widetilde{\N}:=(\sigma(\N'))'\subseteq {\mathcal B}(K)\,. \]
Next, we check that $WW^*\in \widetilde{\N}$\,. Given any $a\in \M$
and $x\in H$, we have by the definition of $W$ that
\[ \langle W^*[a\otimes x]\,, y\rangle=\langle [a\otimes x]\,,
Wy\rangle=\langle [a\otimes x]\,, [1\otimes y]\rangle=\langle
\alpha(a)x\,, y\rangle\,, \quad y\in H\,. \] Hence $W^*[a\otimes
x]=\alpha(a)x$\,, so $WW^*[a\otimes x]=[1\otimes \alpha(a)x]$\,,
which implies that \[ \sigma(n')WW^*[a\otimes x]=[1\otimes
n'\alpha(a)x]\,. \] Moreover, by (\ref{eq7777644443233}),
$WW^*\sigma(n')[a\otimes x]=[1\otimes \alpha(a)n'x]=[1\otimes
n'\alpha(a)x]$\,. We conclude that $WW^*\in \widetilde{\N}$\,. Note
that
\begin{equation}\label{eq4444444446555444}
Wn'=\sigma(n')W\,, \quad n'\in \N'\,,
\end{equation}
since for all $x\in H$\,, $Wn'x=[1\otimes n'x]$\,, while
$\sigma(n')Wx=\sigma(n')[1\otimes x]=[1\otimes n'x]$\,. Taking
adjoints we get
\begin{equation}\label{eq44444444465554443}
n'W^*=W^*\sigma(n')\,, \quad n'\in \N'\,.
\end{equation}
We next prove that $WW^*$ is a properly infinite projection in
$\widetilde{\N}$\,. Define
\[ \rho(x)=WxW^*\,, \quad x\in \N\,. \]
Since $W^*W=1_{\N}$\,, $\rho$ is a (non-unital) $*$-homomorphism of
$\N$ into ${\mathcal B}(K)$\,. Moreover,
\[ WxW^*\sigma(n')=Wxn'=Wn'xW^*=\sigma(n')WxW^*\,, \quad x\in \N\,, n'\in \N'\,. \]
Hence $\rho(x)\in (\sigma(\N'))'=\widetilde{N}$\,, for all $x\in
\N$\,. Since $\rho(1_{\N})=WW^*$\,, we can consider $\rho$ as a
unital $*$-homomorphism of $\N$ into the corner algebra
$(WW^*)\widetilde{N}(WW^*)$\,, and since $\N$ is properly infinite,
it follows that $\rho(\N)$ and $(WW^*)\widetilde{N}(WW^*)$ are also
properly infinite von Neumann algebras. Hence $WW^*$ is a properly
infinite projection in $\widetilde{\N}$\,. Now let $c(WW^*)$ be the
central support of $WW^*$ in $\widetilde{\N}$ and put
$q:=1-c(WW^*)$\,. Then $q\in Z(\widetilde{N})=Z(\sigma(\N'))$ and
therefore $q=\sigma(q_0)$\,, for a projection $q_0\in Z(\N')$\,.
Since $\sigma(q_0)$ and $WW^*$ are orthogonal projections, we have
for all $x\in H$ that
\[ Wq_0x=WW^*[1\otimes q_0x]=WW^*\sigma(q_0)[1\otimes x]=0\,. \]
Therefore $q_0=0$\,, which implies that $c(WW^*)=1_K$\,. By Lemma
\ref{lem82988298}\,, it follows that $WW^*\sim 1_K$ (in
$\widetilde{\N}$)\,. Choose now $U\in \widetilde{\N}$ such that
$U^*U=WW^*$ and $UU^*=1_K$\,. Then $UW\in {\mathcal B}(H, K)$ and
\[ (UW)^*UW=1\,, \quad UW(UW)^*=1_K\,, \]
i.e., $UW$ is a unitary operator from $H$ to $K$\,. Therefore
\[ \pi(a):=(UW)^*\pi_0(a)UW\,, \quad a\in \M \]
defines a normal unital $*$-homomorphism of $\M$ into ${\mathcal
B}(H)$\,. For all $a\in \M$\,,
\[ U^*\pi_0(a)U\in \widetilde{\N}=(\sigma(\N'))'\,. \]
Hence, by (\ref{eq4444444446555444}) and
(\ref{eq44444444465554443})\,, $\pi(a)=W^*(U^*\pi_0(a)U)W\in
(\N')'=\N$\,. Next, set $V:=W^*U^*W\in {\mathcal B}(H)$\,. Since
$U^*\in \widetilde{\N}=(\sigma(\N'))'$\,, using again
(\ref{eq4444444446555444}) and (\ref{eq44444444465554443}) we deduce
that $V\in (\N')'=\N$\,. Moreover, for all $a\in \M$\,,
\[ \alpha(a)=W^*\pi_0(a)W=W^*UW\pi(a)(UW)^*W=V^*\pi(a)V\,. \]
This completes the proof.
\end{proof}

Next we prove the following Wittstock-Haagerup-Paulsen-type theorem
(see \cite{Pau}, Theorem 7.4 and the references given therein):

\begin{theorem}\label{th5}
Let $\M$ and $\N$ be von Neumann algebras with separable preduals
$\M_*\,, \N_*$. Assume, moreover, that $\N$ is properly infinite and
injective. Let $\alpha: \M\rightarrow \N$ be a normal cb-map. Then
there exists a normal $*$-representation $\pi: \M\rightarrow \N$ and
operators $R\,, S\in \N$ such that
\begin{equation}\label{eq20}
\alpha(a)=R\pi(a)S\,, \quad a\in \M \end{equation} and
$\|R\|\|S\|=\|\alpha\|_{\text{cb}}$\,.
\end{theorem}

For the proof we need a few preliminary considerations. Recall that
(see, e.g., \cite{Tk} Vol.I, Theorem 2.14) if $\M$ is a von Neumann
algebra, then any functional $\phi\in \M^*$ has a unique
decomposition into its normal and singular part
\begin{equation}\label{eq17}\phi=\phi_n+\phi_s\,. \end{equation}
Moreover, if $\M$ and $\N$ are von Neumann algebras (not necessarily
with separable preduals), and $T\in {\mathcal B}(\M, \N)$\,, then
there exists a unique decomposition
\begin{equation}\label{eq15} T=T_n+T_s\,,
\end{equation}
where $T_n\,, T_s\in {\mathcal B}(\M, \N)$\,, and such that for any
$\phi\in \N_*$\,, we have
\begin{equation}\label{eq16}
\phi\circ T_n=(\phi\circ T)_n\,, \quad \phi\circ T_s=(\phi\circ
T)_s\,, \end{equation} (cf. \cite{To}, Theorem 1). Further, the
following assertions hold:
\begin{enumerate}
\item [$(a)$] If $T$ is positive, then both $T_n$ and $T_s$ in
the decomposition (\ref{eq15}) are positive.
\item [$(b)$] If $T$ is completely positive, then both $T_n$ and $T_s$
in the decomposition (\ref{eq15}) are completely positive.
\end{enumerate}
Statement $(a)$ follows by uniqueness, since for any positive
functional $\phi\in \M^*$, both maps $\phi_n$ and $\phi_s$ in the
decomposition (\ref{eq17}) are positive. To justify $(b)$\,, let $k$
be a positive integer, and note that for any $\phi\in (M_k(\M))^*$,
we have $\phi=(\phi_{ij})_{i, j=1}^k$ with $\phi_{ij}\in \M^*$\,,
since, algebraically, $(M_k(\M))^*=M_k(\M^*)$. By uniqueness we
deduce that $\phi_n=((\phi_{ij})_n)_{i, j=1}^k$ and, respectively,
$\phi_s=((\phi_{ij})_s)_{i, j=1}^k$.  Thus $(b)$ follows.\\

{\em Proof of Theorem \ref{th5}}. By Definition 1.1 and Theorem 1.6
in \cite{Ha1}, there exist completely positive maps $\beta, \gamma:
\M\rightarrow \N$ such that $\|\beta\|\leq
\|\alpha\|_{\text{cb}}$\,, $\|\gamma\|\leq \|\alpha\|_{\text{cb}}$
and the mapping $\sigma$ defined by
\begin{equation}\label{eq21} \sigma(a):=\left(
\begin{array}
[c]{cc}%
\beta(a) & \alpha^*(a)\\
\alpha(a) & \gamma(a)
\end{array}
\right)\,, \quad a\in \M \end{equation} is a completely positive map
from $\M$ into $M_2(\N)$\,, where
\[ \alpha^*(a):=\alpha(a^*)^*\,, \quad a\in \M\,. \]
Next, note that the maps $\beta$ and $\gamma$ can be chosen to be
normal. For this, exchange (possibly) $\beta$ and $\gamma$ above
with their normal parts $\beta_n$ and $\gamma_n$\,, respectively.
Then, by assertion $(b)$ above one can check that the map $\sigma$
defined by (\ref{eq21}) is still completely positive, and, moreover,
normal. By Theorem \ref{th33}, there exists a normal
$*$-representation ${\pi}': \M\rightarrow M_2(\N)$ and an operator
$V\in M_2(\N)$ such that
\[ \left(
\begin{array}
[c]{cc}%
\beta(a) & \alpha^*(a)\\
\alpha(a) & \gamma(a)
\end{array}
\right)=V^*{\pi}'(a)V\,, \quad a\in \M\,. \] Write now $V=\left(
\begin{array}
[c]{cc}%
V_{11} & V_{12}\\
V_{21} & V_{22}
\end{array}
\right)$\,. It then follows that $\|V\|^2\leq \max\{\|\beta\|\,,
\|\gamma\|\}\leq \|\alpha\|_{\text{cb}}$ and
\[ \alpha(a)=(V_{12}^*\,, V_{22}^*)\,{\pi}'(a)\left(
\begin{array}
[c]{c}%
V_{11} \\
V_{21}
\end{array}
\right)\,, \quad a\in \M\,. \] Since $\N$ is properly infinite,
$\N\cong M_2(\N)$\,. Denote by $1$ the identity of $\N$ and choose
isometries $u_1\,, u_2\in \N$ so that $u_1u_1^*$ and $u_2u_2^*$ are
orthogonal projections with $u_1u_1^*+ u_2u_2^*=1$\,. Define
${\pi}:\M\rightarrow \N$ by \[ {\pi}(a):=(u_1\,,
u_2)\,{\pi}'(a)\left(
\begin{array}
[c]{c}%
u_1^* \\
u_2^*
\end{array}
\right)\in \N\,, \quad  a\in \M\,. \] Then ${\pi}$ is a
$*$-representation because $(u_1\,, u_2)$ is a unitary from $H$ to
$H\oplus H$\,, as verified by the following computations: $(u_1\,,
u_2)\left(
\begin{array}
[c]{c}%
u_1^* \\
u_2^*
\end{array}
\right)=u_1u_1^*+u_2u_2^*=1$\,, respectively, $\left(
\begin{array}
[c]{c}%
u_1^* \\
u_2^*
\end{array}
\right)(u_1\,, u_2)=\left(
\begin{array}
[c]{cc}%
1 & 0\\
0 & 1
\end{array}
\right)$\,. Hence \[ {\pi}'(a)=\left(
\begin{array}
[c]{c}%
u_1^* \\
u_2^*
\end{array}
\right){\pi}(a)(u_1\,, u_2)\,, \quad a\in \M\,. \] Moreover, it is
clear that ${\pi}$ is normal, and that
\begin{equation}\label{eq22} \alpha(a)=(V_{12}^*\,, V_{22}^*)\left(
\begin{array}
[c]{c}%
u_1^* \\
u_2^*
\end{array}
\right){\pi}'(a)(u_1\,, u_2)\left(
\begin{array}
[c]{c}%
V_{11} \\
V_{21}
\end{array}
\right)\,, \quad a\in \M\,. \end{equation} Denote $(V_{12}^*\,,
V_{22}^*)\left(
\begin{array}
[c]{c}%
u_1^* \\
u_2^*
\end{array}
\right)$ by $R$\,, respectively, $(u_1\,, u_2)\left(
\begin{array}
[c]{c}%
V_{11} \\
V_{21}
\end{array}
\right)$ by $S$\,. Then (\ref{eq22}) yields (\ref{eq20}). Moreover,
$\|R\|\|S\|\leq \|V\|^2\leq \|\alpha\|_{\text{cb}}$\,. By
(\ref{eq20}), the reverse inequality $\|\alpha\|_{\text{cb}}\leq
\|R\|\|S\|$ holds, as well, and the proof of theorem \ref{th5} is
complete.\qed

The following result is known as Pelczynski's trick (cf. \cite{Pel};
see also \cite{LT}, page 54):
\begin{lemma}\label{lem1}
Let $X$ and $Y$ be Banach spaces. Suppose that there exist Banach
spaces $V$ and $W$ such that
\begin{enumerate}
\item [$1)$] $X\cong Y\oplus W$
\item [$2)$] $Y\cong X\oplus V$
\item [$3)$] $X\oplus X\cong X$
\item [$4)$] $Y\oplus Y\cong Y$
\end{enumerate}
Then $X$ is isomorphic to $Y$\,.
\end{lemma}

\begin{proof}
For completeness, we include the short proof of this result. We have
\[ X\cong X\oplus X\cong Y\oplus Y\oplus W\oplus W\cong Y\oplus
W\oplus W\cong X\oplus W\,, \] and therefore $X\cong Y\oplus W\cong
(X\oplus V)\oplus W\cong (X\oplus W)\oplus V\cong X\oplus V\cong
Y$\,, as wanted.
\end{proof}

\begin{rem}\label{rem1}\rm
A similar proof with isomorphisms being replaced by complete
isomorphisms shows that Pelczynski's trick holds, more generally, in
the category of operator spaces. Note also that if $X$ and $Y$ are
properly infinite von Neumann algebras, or preduals of properly
infinite von Neumann algebras, then conditions $3)$ and $4)$ above
are automatically satisfied, in the operator space category, as
justified in the proof of Theorem 6.2 in \cite{HRS}.
\end{rem}

\begin{prop}\label{prop1}
Let $\M$ and $\N$ be von Neumann algebras. The following statements
are equivalent:
\begin{enumerate}
\item [$1)$] There exists a cb-embedding $i: \M_*\hookrightarrow \N_*$
such that $i(\M_*)$ is cb-complemented in $\N_*$\,.
\item [$2)$] There exist cb maps $\phi: \M_*\rightarrow \N_*$ and $\psi:
\N_*\rightarrow \M_*$ such that $\psi\circ \phi=\text{Id}_{\M_*}$\,.
\item [$3)$] There exist normal cb-maps $\alpha:\M\rightarrow \N$ and
$\beta:\N\rightarrow \M$ such that $\beta\circ
\alpha=\text{Id}_{\M}$\,.
\end{enumerate}
\end{prop}

\begin{proof}
$1) \Rightarrow 2)$. By hypothesis, there exists a cb-projection
$\rho: \N_*\rightarrow i(\M_*)\subseteq \N_*$\,. Set $\phi:=i$\,,
$\psi:=i^{-1}\circ \rho$\,. Then $\phi, \psi$ are cb-maps and
$\psi\circ \phi=\text{Id}_{\M_*}$\,.

\noindent $2)\Rightarrow 1)$. Set $i:=\phi, \rho:=\phi\circ \psi$\,.
Note that $i=\phi$ is a cb-embedding, since $\phi:\M_*\rightarrow
\phi(\M_*)$ is a bijection with inverse
$\phi^{-1}=\psi\vert_{\phi(\M_*)}$\,. Also, $\rho^2=\phi\circ
(\psi\circ \phi)\circ \psi=\phi\circ \psi=\rho$\,, and moreover,
$\rho(i(\M_*))=\phi(\M_*)$\,. Hence $\rho$ is a cb projection of
$\N_*$ onto $i(\M_*)=\phi(\M_*)$\,.

\noindent $2)\Rightarrow 3)$. Set $\alpha:=\psi^*:\M\rightarrow \N$
and $\beta:=\phi^*: \N\rightarrow \M$\,. Then $\alpha\,, \beta$ are
normal cb-maps with $\beta\circ \alpha=\text{Id}_{\M}$\,. To prove
$3)\Rightarrow 2)$\,, take $\phi:=\beta_*\,, \psi:=\alpha_*$\,. Then
$\phi\,, \psi$ are well-defined cb-maps with $\psi\circ
\phi=\text{Id}_{\M_*}$\,.
\end{proof}

\begin{lemma}\label{lem6}
Let $\M$ be an injective von Neumann algebra with separable predual.
Let $\pi:\M\rightarrow {\mathcal B}(H)$ and $\rho:\M\rightarrow
{\mathcal B}(K)$ be two normal, unital $*$-representations of $\M$
on separable Hilbert spaces $H$ and $K$\,, respectively. Then, for
all $T\in {\mathcal B}(H, K)$\,, there exists $T_0\in
{\overline{\text{conv}\{{\rho(u)}T\pi(u)^*; u\in {\mathcal
U}(\M)\}}}^{w^*},$ where ${\mathcal U}(\M)$ is the set of unitaries
in $\M$\,, such that
\[ T_0=\rho(u)T_0\pi(u)^*\,, \quad u\in {\mathcal U}(\M)\,.
\]
\end{lemma}

\begin{proof}
Let us first recall the following definition due to J. Schwartz (cf.
\cite{Sch}; see also \cite{Sak}). A von Neumann algebra $\N\subseteq
{\mathcal B}(L)$ has property $P$ if ${\overline{\text{conv}\{uTu^*;
u\in {\mathcal U}(\N)\}}}^{w^*}\cap {\N}'\neq \emptyset$\,, for all
$T\in {\mathcal B}(L)$\,, where ${\N}'$ is the commutant of $\N$\,.
It was proved by Sakai that if $\N$ is hyperfinite, then $\N$ has
property $P$ (see \cite{Sak}, Corollary 4.4.19).

Now, let $T\in {\mathcal B}(H, K)$\,. Then, for any $u\in {\mathcal
U}(\M)$\,,
\begin{equation}\label{eq35}
\left(
\begin{array}
[c]{cc}%
\pi(u) & 0\\
0 & \rho(u)
\end{array}
\right)\left(
\begin{array}
[c]{cc}%
0 & 0\\
T & 0
\end{array}
\right)\left(
\begin{array}
[c]{cc}%
\pi(u)^* & 0\\
0 & \rho(u)^*
\end{array}
\right)=\left(
\begin{array}
[c]{cc}%
0 & 0\\
\rho(u)T\pi(u)^* & 0
\end{array}
\right)\,.
\end{equation}
Since $\pi$ and $\rho$ are normal, it follows that the von Neumann
algebra $(\pi\oplus \rho)(\M)$ is injective, and therefore, by the
above discussion, it has property $P$\,. It follows that there
exists an operator \[ T_0\in
{\overline{\text{conv}\{{\rho(u)}T\pi(u)^*; u\in {\mathcal
U}(\M)\}}}^{w^*} \] such that $\left(
\begin{array}
[c]{cc}%
0 & 0\\
T & 0
\end{array}
\right)\in ((\pi\oplus \rho)(\M))'$\,. By applying (\ref{eq35}) to
the operator $T_0\in {\mathcal B}(H, K)$\,, we deduce that
\[ T_0=\rho(u)T_0\pi(u)^*\,, \quad u\in {\mathcal U}(\M)\,,
\]
and the proof is complete.
\end{proof}

\begin{prop}\label{thprop4}
Let $\M$ and $\N$ be properly infinite von Neumann algebras with
separable preduals $\M_*\,, \N_*$. If $\N$ is injective
(=hyperfinite), then the following statements are equivalent:
\begin{enumerate}
\item [$i)$] There exists a cb-embedding $i: \M_*\hookrightarrow \N_*$
such that $i(\M_*)$ is cb-complemented in $\N_*$\,.
\item [$ii)$] There exists a von Neumann algebra embedding $\alpha:
\M\hookrightarrow \N$ and a normal conditional expectation $\beta:
\N\rightarrow \M$ such that $\beta\circ \alpha=\text{Id}_{\M}$\,.
\end{enumerate}
Moreover, if $i)$ holds, then $\M$ is injective, as well.
\end{prop}

\begin{proof} We have to prove that
$i)\Rightarrow ii)$\,, since by Proposition \ref{prop1} we know
already that $ii)\Rightarrow i)$\,. Suppose that $i)$ holds, then by
Proposition \ref{prop1} there exist normal completely bounded maps
$\alpha_1: \M\rightarrow \N$ and $\beta_1: \N\rightarrow \M$ such
that $\beta_1\circ \alpha_1=\text{Id}_{\M}$\,. The first goal is to
make $\alpha_1$ into a $*$-representation. Indeed, by Theorem
\ref{th5} there exists a normal $*$-representation
$\widehat{\pi}:\M\rightarrow \N$ and operators $R, S\in \N$ such
that
\begin{equation}\label{eq567}
\alpha_1(a)=R\widehat{\pi}(a)S\,, \quad a\in \M\,.
\end{equation}
Note that $\alpha_1$ is one-to-one (since $\beta_1\circ
\alpha_1=\text{Id}_{\M}$), and by (\ref{eq567}) this implies that
$\widehat{\pi}$ is one-to-one, too. Hence $\widehat{\pi}$ is a
$*$-isomorphism of $\M$ onto its image $\widehat{\pi}(\M)$\,, the
latter being a von Neumann subalgebra of $\N$\,. Set
\[ \rho(b):=(\widehat{\pi}\circ \beta_1)(RyS)\,, \quad y\in \N\,. \]
Then, for all $b=\widehat{\pi}(a)$\,, where $a\in \M$\,,
$\rho(b)=(\widehat{\pi}\circ
\beta)(R\widehat{\pi}(a)S)=\widehat{\pi}\circ \beta_1\circ
\alpha_1(a)=\widehat{\pi}(a)=b$\,. Thus
$\rho(\widehat{\pi}(\M))=\widehat{\pi}(\M)$\,. Since
$\rho(\N)\subseteq \widehat{\pi}(\M)$\,, we infer that
$\rho(\N)=\widehat{\pi}(\M)$\,, i.e., $\rho$ is a projection of $\N$
onto $\widehat{\pi}(\M)$\,. Hence we have proved that there exist a
normal one-to-one $*$-representation $\widehat{\pi}:\M\rightarrow
\N$ and a normal cb-projection $\rho:\N\rightarrow
\widehat{\pi}(\M)$\,. Let $H$ be a separable Hilbert space with
$\N\subseteq {\mathcal B}(H)$\,. Since $\N$ is injective, there
exists a conditional expectation $E:{\mathcal B}(H)\rightarrow
\N$\,. The composition ${\hat{\pi}}^{-1}\rho E: {\mathcal
B}(H)\rightarrow \M$ is a cb-projection. By a result of Pisier (cf.
Theorem 2.9 in \cite{Pi6}) and Christensen-Sinclair \cite{ChS2}, it
follows that $\M$ is injective, as well.

So far we have reduced the general case to the case when $\M$ is a
von Neumann subalgebra of $\N$\,, $\alpha:\M\rightarrow \N$ is the
inclusion map and $\beta_1:\N\rightarrow \M$ is a normal
cb-projection of $\N$ onto $\M$\,. The next step is to change
$\beta_1$ into a normal conditional expectation. For this, apply now
Theorem \ref{th5} to the map $\beta_1:\N\rightarrow \M$\,, and infer
the existence of a normal $*$-representation $\pi:\N\rightarrow \M$
and operators $R, S\in \M$ such that $\beta_1(a)=R\pi(a)S$\,, for
all $a\in \N$\,. Since $\beta_1\vert_{\M}=\text{Id}_{\M}$ and
$\M\subseteq \N$\,, it follows that
\begin{equation}\label{eq33223322331111}
a=R\pi(a)S\,, \quad a\in \M\,.
\end{equation}
Since $\M$ is injective, we get from Lemma \ref{lem6} that there
exists an operator
\[
S_0\in {\overline{\text{conv}\{{\pi(u)}Su^*; u\in {\mathcal
U}(\M)\}}}^{w^*}\subseteq \M \] such that
\begin{equation}\label{eq332233223311118}
\pi(u)S_0u^*=S_0\,, \quad u\in {\mathcal U}(\M)\,.
\end{equation}
Note that by (\ref{eq33223322331111}), $R\pi(u)Su^*=uu^*=1$\,, for
all $u\in {\mathcal U}(\M)$\,. Hence $RS_0=1$\,. Therefore $S_0$ is
bounded away from 0\,, i.e., $|S_0|:=(S_0^*S_0)^{\frac12}$ is
invertible. Let \[ S_0:=U_0|S_0| \] be the polar decomposition of
$S_0$\,. Since $S_0\in \M$\,, it follows that $U_0\in \M$ and
$|S_0|\in \M$\,, as well. Moreover, $U_0=S_0|S_0|^{-1}$ and
$U_0^*U_0=1$\,. By (\ref{eq332233223311118}), $\pi(u)S_0=S_0u$\,,
for all $u\in {\mathcal U}(\M)$\,. Hence $\pi(a)S_0=S_0a$\,, for all
$a\in \M=\text{Span}({\mathcal U}(\M))$\,. By taking adjoints, it
follows that $S_0^*\pi(a)=aS_0^*$\,, for all $a\in \M$\,. Hence
\[ S_0^*S_0a=S_0^*\pi(a)S_0=aS_0^*S_0\,, \quad a\in \M\,. \]
Therefore $|S_0|=(S_0^*S_0)^{\frac12}\in Z(\M)$ (the center of
$\M$), which implies that
\[
\pi(a)U_0=\pi(a)S_0|S_0|^{-1}=S_0a|S_0|^{-1}=S_0|S_0|^{-1}a=U_0a\,,
\quad a\in \M\,. \] Hence $a=U_0^*\pi(a)U_0$\,, for all $a\in \M$\,,
and therefore $a=U_0^*\pi(a)U_0$\,, for all $a\in \M$\,. From this
we infer that the map $\beta$ defined by
\[ \beta(b):=U_0^*\pi(b)U_0\,, \quad b\in \N \]
is a normal, completely positive map satisfying $\beta(\N)\subseteq
\M$ and ${\beta}\vert_{\M}=\text{Id}_{\M}$\,. Hence $\beta$ is a
normal conditional expectation of $\N$ onto $\M$\,. This completes
the proof of Proposition \ref{thprop4}.
\end{proof}

By Proposition \ref{prop1} and Proposition \ref{thprop4} we obtain
immediately the following
\begin{theorem}\label{th7}
Let $\M$ and $\N$ be properly infinite von Neumann algebras with
separable preduals $\M_*\,, \N_*$. If both $\M$ and $\N$ are
injective, then the following statements are equivalent:
\begin{enumerate}
\item [$1)$] $\M_*\overset{\text{c.b.}}{\cong} \N_*$
\item [$2)$] $\text{Id}_{\M_*}$ admits a cb-factorization through
$\N_*$, and $\text{Id}_{\N_*}$ admits a cb-factorization through
$\M_*$\,.
\item [$3)$] $\text{Id}_{\M}$ admits a cb-factorization through
$\N$ and $\text{Id}_{\N}$ admits a cb-factorization through $\M$\,,
where all four cb-maps involved are normal.
\item [$4)$] There exist von Neumann algebras embeddings $i:
\M\hookrightarrow \N$\,, $j:\N\hookrightarrow \M$ and normal
conditional expectations $E: \N\rightarrow i(\M)$\,,
$F:\M\rightarrow j(\N)$\,.
\end{enumerate}
\end{theorem}

The following result is due to Sakai and Tomiyama (cf. \cite{Sak2}
and \cite{To}, Theorem 3). For convenience, we include a short proof
based on \cite{Sak2}.

\begin{lemma}\label{lem767654}(\,\cite{Sak2}, \cite{To}\,)
Let $\N$ be a  semifinite von Neumann algebra and let $\M\subseteq
\N$ be a von Neumann subalgebra of type III\,. Then there is no
normal conditional expectation from $\N$ onto $\M$\,.
\end{lemma}

\begin{proof}
Suppose by contradiction that there exists a normal conditional
expectation $E:\N\rightarrow \M$\,. Since $\N$ is semifinite, there
exists a net $(e_\lambda)_{\lambda\in \Lambda}$ of finite
dimensional projections in $\N$ converging in strong operator (s.o.)
topology to the identity $1_{\N}$ of $\N$. By normality of $E$\,, it
follows that $E(e_\lambda)\overset{\text{s.o.}}{\longrightarrow}
1$\,, where $1$ is the identity of $\M$\,. Hence, there exists a
finite projection $e\in \N$ with $E(e)\neq 0$\,. Moreover, we can
choose $\varepsilon> 0$ such that \[ p:=1_{[\varepsilon,
\infty)}(E(e))\neq 0\,. \] We show next that $p$ is a finite
projection in $\M$\,. Set $a:=E(e)p+(1-p)$\,. Then $a=a^*$ and
$a\geq \varepsilon 1$\,. In particular, $a$ is invertible. Let
$(x_\alpha)_{\alpha\in A}$ be a bounded net in $p{\M}p$\,, which
converges s.o. to 0\,. Then, as shown by Sakai (see the proof of
Theorem 2.5.6 $(2)$ in \cite{Sak}), the finiteness of $e$ ensures
that the net $(ex_\alpha^*)_{\alpha\in A}$ converges to 0 s.o. Since
$x_\alpha^*\in p{\M}p$\,, for all $\alpha\in A$\,, it follows that
$ax_\alpha^*=E(e)x_\alpha^*=E(ex_\alpha^*)$\,. Thus
$x_\alpha^*=a^{-1}E(ex_\alpha^*)\overset{\text{s.o.}}{\longrightarrow}
0$\,. By Theorem 2.5.6 $(1)$ in \cite{Sak}, it follows that $p{\M}p$
is finite. This implies that $p$ is a finite projection in $\M$,
which contradicts the fact that $\M$ is of type III.
\end{proof}

{\em Proof of Theorem \ref{th1}}. Let $\M$ and $\N$ be hyperfinite
von Neumann algebras with separable preduals, where $\M$ is of type
III and $\N$ is semifinite.

Consider first the case when $\N$ is properly infinite, and assume
by contradiction that $\M_*$ is cb-isomorphic to a cb-complemented
subspace of $\N_*$\,. Then the statement $ii)$ in Proposition
\ref{thprop4} holds. By Lemma \ref{lem767654}, this yields a
contradiction, and therefore the theorem is proved in this case.

For the general case, note that if $\M_*$ is cb-complemented in
$\N_*$\,, then also $(\M\otimes {\mathcal B}(H))_*$ is
cb-complemented in $(\N\otimes {\mathcal B}(H))_*$\,, where
$H=l^2(\mathbb{N})$\,. But the von Neumann algebras $\M\otimes
{\mathcal B}(H)$ and $\N\otimes {\mathcal B}(H)$ are both properly
infinite; moreover, the first is of type III, while the second is
semifinite. So, by the first part of the proof we obtain a
contradiction, and the proof of Theorem \ref{th1} is complete.\qed

\begin{rem}\label{rem444321}\rm
In \cite{ORS}, Oikhberg, Rosenthal and St${\o}$rmer studied {\em
isometric} embeddings of a predual $\M_*$ into a predual $\N_*$\,,
and show that no such embeddings exist if $\M$ is of type III and
$\N$ is semifinite. In this case neither hyperfiniteness, nor
separability assumptions are needed.
\end{rem}

\section{A characterization of type III-factors whose preduals are cb-isomorphic to the predual of the injective type $III_1$-factor}
\setcounter{equation}{0}

Our main result in this section is a characterization of type III-
factors which admit an invariant normal state on their flow of
weights. (See Theorem \ref{th8} below.) Let $\M$ be a type
III-factor with separable predual and let $\phi_0$ be a fixed normal
faithful state on $\M$\,. Consider the crossed-product
$\N:=\M\rtimes_{\sigma^{\phi_0}} \mathbb{R}$, and let
$(\widetilde{\theta}_s)_{s\in \mathbb{R}}$ of $\mathbb{R}$ be the
dual action of $\sigma^{\psi_0}$ on $\N$\,, i.e.,
\begin{eqnarray}
\widetilde{\theta}_s(\pi(x))&=&\pi(x)\,, \quad x\in \M\label{eq333334333339}\\
\widetilde{\theta}_s(\lambda(t))&=&e^{ist} \lambda(t)\,, \quad t\in
\mathbb{R}\,,\label{eq3333343333392}
\end{eqnarray}
for all $s\in \mathbb{R}$\,, where $\pi(\M)$ and $(\lambda(t))_{t\in
\mathbb{R}}$ are the generators of the crossed product $\N$, as
explained in the Introduction. Then the flow of weights of $\M$ is
the pair $(Z(\N)\,, (\theta_s)_{s\in \mathbb{R}})$\,, where $Z(\N)$
is the center of $\N$ and $\theta_s$ is the restriction of
$\widetilde{\theta_s}$ to the center $Z(\N)$\,, for all $s\in
\mathbb{R}$\,. Recall that a state $\psi$ on $Z(\N)$ is called {\em
an invariant state on the flow of weights for} $\M$, if $\psi$ is
invariant under the dual action $(\theta_s)_{s\in \mathbb{R}}$\,,
i.e.,
\[ \psi\circ \theta_s=\psi\,, \quad s\in \mathbb{R}\,. \]
Next, let $h$ be the unique positive self-adjoint operator
affiliated with $\N$, for which $h^{it}=\lambda(t)$\,, for all $t\in
\mathbb{R}$\,, and let $\widetilde{\phi}_0$ be the dual weight of
$\phi_0$ in the sense of Takesaki \cite{Ta2}. Then (cf. \cite{Ta2}),
$\widetilde{\phi}_0=\tau(h\,\cdot \,)$\,, for a unique normal,
faithful, semifinite trace $\tau$ on $\N$\,, which satisfies
\[ \tau\circ \widetilde{\theta}_s=e^{-s}\tau\,, \quad s\in
\mathbb{R}\,.
\] By \cite{HS}, we can define a map $\phi\mapsto \widehat{\phi}$ of
$S_{\text{nor}}(\M)$\,, the set of normal states on $\M$\,, into the
set $S_{\text{nor}}(Z(\N))$ of normal states on $Z(\N)$ in the
following way: For $\phi\in S_{\text{nor}}(\M)$\,, let $h_\phi$
denote the unique positive unbounded operator affiliated with $\N$,
such that \begin{equation}\label{eq33445544}
\widetilde{\phi}=\tau(h_\phi \,\cdot \,)\,, \end{equation} where
$\widetilde{\phi}$ is the dual weight of $\phi$, and set
\begin{equation}\label{eq998899789}
 e_\phi:=1_{(1,\infty)}(h_\phi)\in \N\,.
\end{equation}
Then $\tau(e_\phi)=1$\,, and therefore
\[ \widehat{\phi}(z)=\tau(e_\phi z)\,, \quad z\in Z(\N) \]
defines a normal state on $Z(\N)$ (cf. \cite{HS}, Def. 3.2). By the
proof of Lemma 3.4 in \cite{HS}\,,
\begin{equation}\label{554555455543}
\widetilde{\theta}_s(e_\phi)=1_{[e^s, \infty)}(h_\phi)\,, \quad s\in
\mathbb{R}\,.
\end{equation}

The Main Theorem in \cite{HS} states that if $\M$ is a properly
infinite von Neumann algebra, then $\phi\mapsto \widehat{\phi}$ maps
$S_{\text{nor}}(\M)$ onto the set $\{\omega\in S_{\text{nor}}(\N);
\, \omega\circ \theta_s\geq e^{-s}\omega\,, \forall s\geq 0\}$\,.
Moreover, if $\phi, \psi\in S_{\text{nor}}(\M)$\,, then
\begin{equation}\label{eq565441}
\|\widehat{\phi}-\widehat{\psi}\|=\inf\limits_{u\in {\mathcal
U}(\M)}\|u\phi u^*-\psi\|\,.
\end{equation}

The right hand-side in (\ref{eq565441}) is by definition (cf.
\cite{CS}) equal to the distance $d([\phi]\,, [\psi])$\,, where
$[\phi]$ denotes the norm-closure of the orbit of $\phi$ under the
action of inner $*$-automorphisms, $\text{Int}(\M)$, by $\phi\mapsto
u\phi u^*$\,. The customary notation $u\phi u^*(x)=\phi(u^*xu)$\,,
for all $x\in \M$ is being used. We write $\phi\sim \psi$ if
$d([\phi]\,, [\psi])=0$\,. We now consider the following more
general equivalence relation:

\begin{defi}\label{def88}\rm Let $m\,, n$ be positive integers. A
$*$-isomorphism $\alpha:M_m(\M)\rightarrow M_n(\M)$ is called {\em
inner} if there exists $u\in M_{n, m}(\M)$ satisfying
$u^*u=1_{M_m(\M)}$ and $uu^*=1_{M_n(\M)}$, such that
\[ \alpha(x)=uxu^*\,, \quad x\in M_m(\M)\,. \]
Let $\text{Int}(M_m(\M), M_n(\M))$ denote the set of all inner
$*$-isomorphisms $\alpha:M_m(\M)\rightarrow M_n(\M)$\,.

If $\phi\in S_{\text{nor}}(M_m(\M))$ and $\psi\in
S_{\text{nor}}(M_n(\M))$\,, then we write $\phi\sim \psi$ if
\[ \inf\{
\|\phi-\psi\circ \alpha\|\,; \alpha\in \text{Int}(M_m(\M),
M_n(\M))\}=0\,. \] Clearly $\sim$ is an equivalence relation on
$\bigcup_{k=1}^\infty S_{\text{nor}}(M_k(\M))$\,.
\end{defi}

\begin{prop}\label{prop99}
Let $\M$ be a properly infinite von Neumann algebra with separable
predual. Given $\phi\in S_{\text{nor}}(\M)$\,, then $\widehat{\phi}$
is $\theta$-invariant if and only if for all $n\geq 1$\,,
\[ \phi\sim \phi_n:={\frac1n}\phi \otimes \text{Tr}_n\,, \]
where $\text{Tr}_n$ denotes the non-normalized trace on
$M_n(\mathbb{C})$\,.
\end{prop}

\begin{proof}
Let $\phi\in S_{\text{nor}}(\M)$\,, and let $n\in \mathbb{N}$\,.
Since $\M$ is properly infinite, we can choose isometries $v_1\,,
\ldots \,, v_n\in \M$ with orthogonal ranges such that $\sum_{i=1}^n
v_iv_i^*=1_{\M}$\,. Then $v:=(v_1\,, \ldots\,, v_n)$ is a unitary in
$M_{1, n}(\M)$\,. We will identify $\M$ with the subalgebra
$\pi(\M)$ of $\N=\M\rtimes_{\sigma^{\phi_0}} \mathbb{R}$\,. Define
now maps $\alpha: \M\rightarrow M_n(\M)$ and
$\widetilde{\alpha}:\N\rightarrow M_n(\N)$ by \begin{eqnarray*}
\alpha(x)&=& v^*xv\,, \quad x\in \M\,,\\
\widetilde{\alpha}(y)&=& v^*yv\,, \quad y\in \N\,,
\end{eqnarray*}
and put
\begin{equation}\label{eq6665556664321}
\psi:=\phi_n\circ \alpha=\left(\phi\otimes
{\frac1n}\text{Tr}_n\right)\circ \alpha\,.
\end{equation}
Then $\psi\sim \phi_n$\,. We claim that the dual weights
$\widetilde{\phi}$ and $\widetilde{\psi}$ on $\N$ of $\phi$ and
$\psi$\,, respectively, satisfy
\begin{equation}\label{eq66655566643218}
\widetilde{\psi}=\left(\widetilde{\phi}\otimes
{\frac1n}\text{Tr}_n\right)\circ \widetilde{\alpha}\,.
\end{equation}
For $x\in \M_+$\,,
\begin{equation}\label{eq666555666432187}
\psi(x)=\left(\phi\circ
{\frac1n}\text{Tr}_n\right)\left((v_i^*xv_j)_{i,
j=1}^n\right)={\frac1n}\sum\limits_{i=1}^n \phi(v_i^*xv_i)\,.
\end{equation}
It can be easily shown that formula (\ref{eq666555666432187}) holds
for $x\in \widehat{\M}_+$\,, the extended positive part of $\M$\,,
as defined in \cite{Ha9}, Sect. 1. Then by \cite{Ha8}, Theorem 1.1,
the dual weights $\widetilde{\phi}$\,, $\widetilde{\psi}$ of $\phi$
and $\psi$ are given by
\begin{eqnarray*}
\widetilde{\phi}(y)&=&\phi\left(\int_\mathbb{R}
\widetilde{\theta}_s(y)\,
ds\right)\,, \quad y\in \N_{+}\,,\\
\widetilde{\psi}(y)&=&\psi\left(\int_\mathbb{R}
\widetilde{\theta}_s(y)\, ds\right)\,, \quad y\in \N_{+}\,,
\end{eqnarray*}
where $\int_\mathbb{R} \widetilde{\theta}_s(y)\, ds$ is an element
of the extended positive part of $\pi(\M)=\M$\,. Note that for all
$s\in \mathbb{R}$ and all $1\leq i\leq n$\,, we have
$\widetilde{\theta}_s(v_i)=v_i$\,, since $v_i\in \M$\,. Hence, by
(\ref{eq666555666432187})\,,
\begin{equation*}
\widetilde{\psi}(y)={\frac1n}\sum\limits_{i=1}^n\widetilde{\phi}(v_i^*yv_i)=\left(\widetilde{\phi}\otimes
{\frac1n}\text{Tr}_n\right)\circ \widetilde{\alpha}(y)\,, \quad y\in
\N_{+}\,,
\end{equation*}
which proves (\ref{eq66655566643218})\,. Next, observe that for
$y\in \N_{+}$\,,
\begin{equation*}
\sum\limits_{i=1}^n \tau(v_i^*yv_i)=\sum\limits_{i=1}^n
\tau(y^{1/2}v_iv_i^*y^{1/2})=\tau(y)\,,
\end{equation*}
because $\sum_{i=1}^n v_iv_i^*=1_{\M}$\,. Hence
\begin{equation}\label{eq6665556664321876}
\tau=\left(\tau\otimes \text{Tr}_n\right)\circ \widetilde{\alpha}\,.
\end{equation}
By (\ref{eq66655566643218}) and (\ref{eq6665556664321876}) we then
have
\begin{equation}\label{eq66655566643218765}
h_{\psi}=\frac{d\widetilde{\psi}}{d\tau}={\widetilde{\alpha}}^{-1}\left(\frac{d\left(\widetilde{\phi}\otimes
{\frac1n}\text{Tr}_n\right)}{d(\tau\otimes \text{Tr}_n)}\right)=
{\frac1n}v(h_\phi\otimes I_n)v^*\,,
\end{equation}
where $I_n\in M_n(\mathbb{C})$ denotes the identity $n\times n$
matrix. Hence by (\ref{554555455543})\,,
\begin{eqnarray*}
e_{\psi}=1_{(1, \infty)}(h_{\psi})=1_{(n, \infty)}(v(h_\phi\otimes
I_n)v^*)
&=& v(1_{(n, \infty)}(h_\phi)\otimes I_n)v^*\\
&=& v(\widetilde{\theta}_{\log n}(e_\phi)\otimes I_n)v^*\,.
\end{eqnarray*}
Therefore, for all $z\in Z(\N)$\,,
\begin{eqnarray*}
\widehat{\psi}(z)=\tau(e_{\psi} z)= \tau(v(\widetilde{\theta}_{\log
n}(e_\phi)\otimes I_n)v^*z)&=& \sum\limits_{i=1}^n
\tau(v_i\widetilde{\theta}_{\log n}
(e_\phi)v_i^* z)\\
&=&\sum\limits_{i=1}^n \tau(v_i^*v_i \widetilde{\theta}_{\log
n}(e_\phi) z)\\
&=& n\tau(\widetilde{\theta}_{\log n}(e_\phi) z)\\&=& n(\tau\circ
\widetilde{\theta}_{\log n})(e_\phi {\theta}_{-{\log n}}(z))\,.
\end{eqnarray*}
Since $\tau\circ \widetilde{\theta}_{\log n}=e^{-{\log
n}}\tau={\frac1n}\tau$\,, it follows that
\begin{equation}\label{eq666555666432187654}
\widehat{\psi}(z)=\tau(e_\phi {\theta}_{-{\log
n}}(z))=(\widehat{\phi}\circ {\theta}_{-{\log n}})(z)\,, \quad z\in
Z(\N)\,.
\end{equation}
By the Main Theorem in \cite{HS} (cf. (\ref{eq565441}) above), we
now deduce that
\[ \psi\sim \phi \,\,\, \Leftrightarrow \,\,\, \widehat{\psi}=\widehat{\phi}
\,\,\, \Leftrightarrow \,\,\, \widehat{\phi}=\widehat{\phi}\circ
\theta_{-{\log n}}\,. \] Since $\psi\sim \phi_n$\,, we get by
transitivity that
\begin{equation*}
\phi\sim \phi_n \,\,\, \Leftrightarrow \,\,\,
\widehat{\phi}=\widehat{\phi}\circ \theta_{-{\log n}}\,.
\end{equation*}
This holds for every $n\in \mathbb{N}$\,. Hence, if $\widehat{\phi}$
is $\theta$-invariant, then $\phi\sim \phi_n$\,, for all $n\in
\mathbb{N}$\,. Conversely, if $\phi\sim \phi_n$ for all $n\in
\mathbb{N}$\,, then
\[ \widehat{\phi}\circ \theta_{-{\log n}}=\widehat{\phi}\,, \quad
 n\in \mathbb{N}\,. \]
 Note that the group generated by $\{-{\log n}; n\geq 1\}$ in
$(\mathbb{R}\,, +)$ is $\log (\mathbb{Q}^+)$\,, which is dense in
$(\mathbb{R}\,, +)$\,. By the continuity of the semigroup
$(\theta_s)_{s\in \mathbb{R}}$\,, we conclude that
$\widehat{\phi}\circ \theta_s=\widehat{\phi}$\,, for all $s\in
\mathbb{R}$\,, i.e., $\widehat{\phi}$ is $\theta$-invariant. The
proof is now complete.
\end{proof}

\begin{cor}\label{cor999}
Let $\M$ be a properly infinite von Neumann algebra with separable
predual. If $\phi\in S_{\text{nor}}(\M)$ such that $\phi\sim
{\frac1n}\phi \otimes \text{Tr}_n$ for all positive integers $n$\,,
then
\begin{equation}\label{eq34342234}\phi\sim \left(
\begin{array}
[c]{cc}%
a\phi & 0\\
0 & (1-a)\phi
\end{array}
\right)\,, \quad 0< a< 1\,.
\end{equation}
\end{cor}

\begin{proof}
Assume first that $a\in \mathbb{Q}$\,. Then $a=\frac{p}{q}$\,, where
$p, q$ are positive integers with $p< q$\,. By hypothesis, we then
have
\[ \left(
\begin{array}
[c]{cc}%
a\phi & 0\\
0 & (1-a)\phi
\end{array}
\right)={\frac1q}\left(
\begin{array}
[c]{cc}%
p\phi & 0\\
0 & (q-p)\phi
\end{array}
\right)\sim {\frac1q}\left(
\begin{array}
[c]{cc}%
\phi\otimes \text{Tr}_p & 0_{p, q-p}\\
0_{q-p, p} & \phi\otimes \text{Tr}_{q-p}
\end{array}
\right)\sim \phi\,, \] where $0_{p, q-p}\in M_{p, q-p}(\mathbb{C})$
and $0_{q-p, p}\in M_{q-p, p}(\mathbb{C})$ denote matrices of
corresponding sizes with all entries equal to zero. By approximation
we obtain (\ref{eq34342234}) for all $0< a< 1$\,.
\end{proof}

\begin{lemma}\label{lem555543212}
Let $\M_1$ and $\M_2$ be factors with separable preduals, and assume
that $\M_2$ is of type III. Then the following two conditions are
equivalent:
\begin{enumerate}
\item [$(a)$] There exists a von Neumann algebra embedding
$i:\M_1\hookrightarrow \M_2$ and a normal conditional expectation
$E:\M_2\rightarrow \M_1$\,.
\item [$(b)$] There exists a von Neumann algebra embedding
$i:\M_1\hookrightarrow \M_2$ and a normal \underline{faithful}
conditional expectation $E:\M_2\rightarrow \M_1$\,.
\end{enumerate}
\end{lemma}

\begin{proof}
The implication $(b)\Rightarrow (a)$ is trivial. We prove that
$(a)\Rightarrow (b)$\,. For this, we can assume that $\M_1\subseteq
\M_2$\,, for which there exists a normal conditional expectation
$E:\M_2\rightarrow \M_1$\,. Since $E$ is normal, it has a support
projection $p:=\text{supp}(E)$\,, $p$ is the smallest projection in
$\M_2$ for which $E(1-p)=0$\,. Moreover, $p\in \M_1^{\prime} \cap
\M_2$ by the bimodule property of $E$\,. Hence the map $\phi$
defined by
\[ \phi(a):=ap=pa\,, \quad a\in \M_1 \]
is a normal $*$-homomorphism of $\M_1$ onto the von Neumann algebra
$p{\M_1}\subseteq {\mathcal B}(pH)$\,, where $\M_1\subseteq
{\mathcal B}(H)$\,. Since $\M_1$ is a factor, the map $\phi$ is
one-to-one, and hence $\M_1\simeq p{\M_1}$\,. Since $\M_2$ is a type
III-factor with separable predual and $p\neq 0$\,, we deduce that
$p\sim 1_{\M_2}$\,, i.e., $p=uu^*$ for an isometry $u\in \M_2$\,.
Hence $\M_2\simeq p{\M_2}p\subseteq {\mathcal B}(pH)$\,. Therefore,
in order to prove $(b)$ it suffices to construct a normal faithful
conditional expectation $E':p{\M_1}p\rightarrow p{\M_2}$\,. Set
\[ E'(y):=pE(y)=E(y)p\,, \quad y\in p{\M_1}p\,. \]
Then $E'$ is positive, unital and normal. For $a, b\in \M_1$ and
$x\in \M_2$\,, we have
\begin{eqnarray*}
E'((pa)(pxp)(bp))&=&pE((ap)(pxp)(pb))p\\&=&pE(a(pxp)b)p\\&=&paE(pxp)bp\,.
\end{eqnarray*}
Hence $E'$ is a normal conditional expectation of $p{\M_2}p$ onto
$p{\M_1}={\M_1}p$\,. Moreover, if $y\in p{\M_2}p$\,, $y\geq 0$ and
$E'(y)=0$\,, then $E(y)=\phi^{-1}(E'(y))=0$\,, and hence the support
projection of $y$ is less than $1-p$\,. It follows that $y=0$\,, and
we have shown that $E'$ is faithful.
\end{proof}

\begin{theorem}\label{th8}
Let $\M$ be a type III-factor with separable predual. The following
statements are equivalent:
\begin{enumerate}
\item [$(1)$] There exists a von Neumann algebra embedding of $R_{\infty}$ into $\M$ with a normal conditional expectation $E:\M\rightarrow
R_\infty$\,.
\item [$(2)$] For every $0< \lambda< 1$\,, there exists a von Neumann algebra embedding of $R_{\lambda}$ into $\M$ with a normal conditional expectation $E:\M\rightarrow
R_{\lambda}$\,.
\item [$(3)$] There exists an invariant normal state on
the flow of weights $(Z(\N), (\theta_s)_{s\in \mathbb{R}})$ for
$\M$\,.
\end{enumerate}
\end{theorem}

\begin{proof}
It was shown by Haagerup, Rosenthal and Sukochev (see \cite{HRS},
Theorem 6.2 $(a)$) that $(R_\lambda)_*$ is completely isomorphic to
$(R_\infty)_*$\,, for all $0< \lambda\leq 1$\,. Then an application
of Theorem \ref{th7} yields the equivalence of statements $(1)$ and
$(2)$\,.

We now prove the implication $(1)\Rightarrow (3)$\,. Following an
argument from \cite{HS}, we will show a (slightly) more general
result. Namely, let $\M_1$ and $\M_2$ be factors of type III with
separable preduals such that there is a von Neumann algebra
embedding of $\M_1$ into $\M_2$ with a normal conditional
expectation $E: \M_2\rightarrow \M_1$\,. We show that if $\M_1$ is a
type III$_1$-factor, then under the above condition there exists an
invariant normal state on the flow of weights for $\M_2$\,. By Lemma
\ref{lem555543212} we can assume that $E$ is faithful. Let $\phi_1$
be a normal, faithful state on $\M_1$ and set $\phi_2:=\phi_1\circ
E$\,. Then $\phi_2$ is a normal faithful state on $\M_2$, and by
\cite{Ta1} it follows that
\begin{equation}\label{eq4567} {\sigma_t^{\phi_2}}_{|_{\M_1}}=\sigma_t^{\phi_1}\,, \quad
t\in \mathbb{R}\,. \end{equation} By (\ref{eq4567}) we obtain an
embedding of $\N_1:={\M_1}\rtimes_{\sigma^{\phi_1}} \mathbb{R}$ into
$\N_2:={\M_2}\rtimes_{\sigma^{\phi_2}} \mathbb{R}$, with a normal
faithful conditional expectation $\widetilde{E}: \N_2\rightarrow
\N_1$\,. Moreover, $\widetilde{E}\vert_{\pi_2(\M_2)}=E$ and
\begin{equation}\label{eq445532211}
\widetilde{E}(\lambda_2(t))=\lambda_1(t)\,, \quad t\in \mathbb{R}\,,
\end{equation}
where $\pi_i(\M_i)$ and $(\lambda_i(t))_{t\in \mathbb{R}}$ are the
generators of $\N_i$\,, $i=1, 2$, as explained before. Further, let
$\theta^{(1)}:=(\theta^{(1)}_s)_{s\in \mathbb{R}}$ and
$\theta^{(2)}:=(\theta^{(2)}_s)_{s\in \mathbb{R}}$ denote the dual
action of $\mathbb{R}$ on $\N_1$ and $\N_2$\,, respectively. By
(\ref{eq333334333339}) and (\ref{eq3333343333392}), $\theta^{(2)}$
extends $\theta^{(1)}$ and the canonical trace $\tau_2$ on $\N_2$
extends the canonical trace $\tau_1$ on $\N_1$\,. Also,
\begin{equation}\label{eq4568}
\widetilde{E}(Z(\N_2))\subseteq Z(\N_1)\,,
\end{equation}
which can be justified as follows. Given $x\in Z(\N_2)$, then
$xy=yx$ for all $y\in \N_1$\,, and therefore by the bimodule
property of conditional expectations,
$y\widetilde{E}(x)=\widetilde{E}(yx)=\widetilde{E}(xy)=\widetilde{E}(x)y$\,.

Since $\M_1$ is a type III$_1$- factor, then $Z(\N_1)=\mathbb{C}
1_{\N_1}$ (cf. {\bf 1.6} in the Introduction)\,, where $1_{\N_1}$ is
the identity of $\N_1$\,. By (\ref{eq4568}) we infer that for all
$x\in Z(\N_2)$ there exists $\phi(x)\in \mathbb{C}$ such that
$\widetilde{E}(x)=\phi(x)\cdot 1_{\N_1}$\,. It is then easily seen
that the correspondence $x\in Z(\N_2)\mapsto \phi(x)\in \mathbb{C}$
defines a normal $\theta^{(2)}$-invariant state on $Z(\N_2)$\,. Thus
assertion $(3)$ is proved.

Now we show the remaining implication $(3)\Rightarrow (1)$\,.
Suppose that $(Z(\N), (\theta_s)_{s\in \mathbb{R}})$ has an
invariant normal state $\omega$\,. Note that, since $\M$ is a
factor, the faithfulness of $\omega$ is automatic. Indeed, we have
that $\theta_s(\text{supp}(\omega))=\text{supp}(\omega)$\,, for all
$s\in \mathbb{R}$\,, where $\text{supp}(\omega)$ denotes the support
projection of $\omega$\,. Since $\M$ is a factor, $(\theta_s)_{s\in
\mathbb{R}}$ is ergodic, which implies that $\omega$ has full
support, and hence $\omega$ is faithful. Clearly,
\[ \omega\circ \theta_s\geq e^{-s}\omega\,, \quad s\geq 0\,.
\]
Hence, by the Main Theorem in \cite{HS}, there exists $\phi\in
S_{\text{nor}}(\M)$ such that $\widehat{\phi}=\omega$\,. By
Proposition \ref{prop99} and Corollary \ref{cor999}, we infer that
for all $0< a< 1$\,,
\begin{equation}\label{eq56564456}
\phi\sim \left(
\begin{array}
[c]{cc}%
a\phi & 0\\
0 & (1-a)\phi
\end{array}
\right)\,.
\end{equation}
Let $0< a_1, a_2< \frac12$ be chosen such that, if
$\lambda_i:=\frac{a_i}{1-a_i}$\,, $i=1, 2$\,, then $\frac{\log
\lambda_1}{\log \lambda_2}\notin \mathbb{Q}$. Further, for all
integers $k\geq 0$ set $a_{2k+1}:=a_1$\,, $a_{2k+2}:=a_2$ and define
\[ \phi_n:=\text{Tr}_2\left(\left(
\begin{array}
[c]{cc}%
a_n & 0\\
0 & (1-a_n)
\end{array}
\right)\,\cdot\,\right)\,, \quad n\geq 1\,. \] Then,
${R_\infty}\cong {R_{\lambda_1}}\otimes {R_{\lambda_2}}\cong
\otimes_{n=1}^\infty (M_2(\mathbb{C}), \phi_n)$ (see {\bf 1.9.} in
the Introduction).

Now let $\varepsilon> 0$\,. By (\ref{eq56564456}), there exists
$\alpha\in \text{Int}(\M, M_2(\M))$ such that $\|\phi-(\phi_1\otimes
\phi)\circ \alpha\|<\frac{\varepsilon}{2}$\,. Hence we can write
$\M=M_2(\mathbb{C})\bar{\otimes} Q_1$\,, and choose a normal
faithful state $\psi_1$ on $Q_1$ so that $(Q_1, \psi_1)\cong (\M,
\phi)$ and
\[ \|\phi-\phi_1\otimes \psi_1\|< \frac{\varepsilon}{2}\,. \]
Similarly, write $Q_1=M_2(\mathbb{C})\bar{\otimes} Q_2$\,, with a
normal faithful state $\psi_2$ on $Q_2$ such that $(Q_2,
\psi_2)\cong (Q_1, \psi_1)$ and $\|\psi_1-\phi_2\otimes \psi_2\|<
{\varepsilon}/{4}$\,. By continuing in this way, we obtain for all
$k\in \mathbb{N}$ a decomposition \[ \M=\left(\otimes_{j=1}^k
M_2(\mathbb{C})\right)\bar{\otimes} Q_k \] and a normal faithful
state $\psi_k$ on $Q_k$ such that $(Q_k, \psi_k)\cong (\M, \phi)$
and, moreover,
\begin{equation}\label{eq565677676}
\|\psi_k-\phi_{k+1}\otimes \psi_{k+1}\|< \frac{\varepsilon}{2^{k+1}}
\,.\end{equation}
Set $\chi_k:=(\phi_1\otimes \ldots \otimes
\phi_k)\otimes \psi_k$\,, for all $k\geq 1$\,. By
(\ref{eq565677676}), we infer that
\begin{eqnarray*}
\|\chi_{k+1}-\chi_k\|&=&\|(\phi_1\otimes \ldots \phi_k)\otimes
(\phi_{k+1}\otimes \psi_{k+1}-\psi_k)\|\\
&=&\|\phi_{k+1}\otimes \psi_{k+1}-\psi_k\|\\&< &\frac
{\varepsilon}{2^{k+1}}\,.
\end{eqnarray*}
This shows that the sequence $(\chi_k)_{k\geq 1}$ converges in
$S_{\text{nor}}(\M)$\,. Set $\chi:=\lim\limits_{k\rightarrow \infty}
\chi_k \in S_{\text{nor}}(\M)$\,. We next show that
$e:=\text{supp}(\chi)$ commutes with $P_k:=\otimes_{j=1}^k
M_2(\mathbb{C})$\,, for all positive integers $k$\,.

Let now $k\in \mathbb{N}$ be fixed. Then, for all $n> k$ we have
\[ \chi_n=(\phi_1\otimes \ldots \otimes \phi_k)\otimes \omega_{k, n}\,, \]
where $\omega_{k, n}:=(\phi_{k+1}\otimes \ldots \otimes
\phi_n)\otimes \psi_n\in S_{\text{nor}}(Q_k)$\,. Further, note that
for all $n, m> k$\,,
\[ \|\chi_n-\chi_m\|=\|\omega_{k, n}-\omega_{k, m}\|\,. \]
Since $\chi:=\lim\limits_{k\rightarrow \infty} \chi_k$\,, we deduce
that $\eta_k:=\lim\limits_{n\rightarrow \infty} \omega_{k, n}$
exists, $\eta_k\in S_{\text{nor}}(Q_k)$ and \[ \chi=(\phi_1\otimes
\ldots \otimes \phi_k)\otimes \eta_k\,. \] Then, since $\phi_1\,,
\ldots ,\phi_k$ are faithful states, it follows that
\[ e:=\text{supp}(\chi)=\text{supp}(\phi_1\otimes \ldots \otimes
\phi_k)\otimes \text{supp}(\eta_k)=1\otimes \text{supp}(\eta_k)\in
1\bar{\otimes} Q_k=P_k^{\prime}\cap \M\,, \] wherein we have used
the fact that $P_k$ is a factor. The claim that $e$ commutes with
$P_k$ is now proved.

Next set $\chi^0:=\chi_{\vert_{e\M e}}$ and define $P_k^0:=eP_k\cong
P_k$\,, $Q_k^0:=eQ_ke$ and $\eta_k^0:={\eta_k}_{\vert_{eQ_k e}}$\,,
for all $k\geq 1$\,. We now obtain a sequence of von Neumann algebra
inclusions $P_1^0\subseteq {P_2^0}\subseteq \ldots \subseteq {e\M
e}$\,. Moreover, for each $k\in \mathbb{N}$\,,
\begin{equation}\label{eq777777766666532313}
e\M e=P_k^0\bar{\otimes} Q_k^0 = (\otimes_{i=1}^k
M_2(\mathbb{C}))\bar{\otimes} Q_k^0\,,
\end{equation}
and with respect to this tensor product decomposition,
\begin{equation}\label{eq7777777666665323134}
\chi^0=(\phi_1\otimes \ldots \otimes \phi_k)\otimes \eta_k^0\,.
\end{equation}
Note that $\chi^0:=\chi_{\vert_{e\M e}}$ is a faithful state on $e\M
e$\,. Let $G$ denote the s.o.t.-closure of $\cup_{k=1}^\infty P_k^0$
in $e\M e$\,. It now follows from (\ref{eq777777766666532313}) and
(\ref{eq7777777666665323134}) that
\[ \left(G, \chi_{\vert_G}^0\right)\cong \otimes_{i=1}^\infty (\M_2(\mathbb{C}),
\phi_i)\cong R_\infty\,,
\]
and for all $k\in \mathbb{N}$ and all $t\in \mathbb{R}$\,,
\[ \left(\sigma_t^{\phi_1\otimes \ldots \otimes
\phi_k}\right)\otimes \sigma_t^{\eta_k^0}=\sigma_t^{\chi^0}\,,
\] which implies that
$(\sigma_t^{\chi^0})_{t\in \mathbb{R}}$ leaves $P_k^0$ globally
invariant for all $k\in \mathbb{N}$\,. It follows that
$(\sigma_t^{\chi^0})_{t\in \mathbb{R}}$ leaves $G$ globally
invariant. By \cite{Ta1} we deduce that there exists a unique
faithful normal conditional expectation $E$ of $e\M e$ onto $G$ such
that $\chi^0\circ E=\chi^0$\,. Since $\M$ is a type III-factor with
separable predual, all non-zero projections in $\M$ are equivalent.
Hence $\M\cong e\M e$\,. Moreover, $G\cong R_\infty$\,. This
completes the proof of the implication $(3)\Rightarrow (1)$\,.
\end{proof}

\begin{rem}\label{rem85645}\rm
In the special case when $\M$ is hyperfinite, a more elementary
proof of the implication $(3)\Rightarrow (1)$ in Theorem \ref{th8}
can be obtained. The proof below was suggested to us by Georges
Skandalis.

Assume that $\M$ is a hyperfinite factor of type III (with separable
predual), such that there exists an invariant normal state $\phi$ on
the flow of weights $(A, \theta^{(1)})$ for $\M$\,. Let
$\N_0:=R_\infty\rtimes_{\sigma^{\omega}} \mathbb{R}$\,, where
$\omega$ is a fixed normal faithful state on $R_\infty$\,. Then
$\N_0$ is the hyperfinite II$_\infty$-factor, and the dual action
$\theta^{(0)}:=\widetilde{\sigma}^\omega$ of $\sigma^\omega$
satisfies $\tau_0\circ \theta^{(0)}=e^{-s}\tau_0$\,, for all $s\in
\mathbb{R}$\,, where $\tau_0$ is a normal faithful trace on
$\N_0$\,. Set $\N:=\N_0\bar{\otimes} A$\,,
$\theta:=\theta^{(0)}\otimes \,\theta^{(1)}$ and
$\tau:=\tau_0\otimes \,\omega$\,. Then $\tau$ is a normal faithful
trace on $\N$ satisfying
\begin{equation}\label{eq22332222243}
\tau\circ \theta_s=e^{-s}\tau\,, \quad s\in \mathbb{R}\,.
\end{equation}
By crossed product theory (cf. \cite{Tk}, Vol. II, Theorem X.2.3
(i)),
\[ R_\infty = \N_0^{\theta^{(0)}}:=\{x\in \N_0;
\,\theta_s^{(0)}(x)=x\,, \text{for all} \,\,s\in \mathbb{R}\}\,. \]
Hence the fixed point algebra $\N^{\theta}$ for the action of
$\theta$ on $\N$ satisfies $R_\infty \bar{\otimes} 1\subseteq
\N^{\theta}\subseteq \N$\,. Put $E:=\text{Id}_{\N_0}\otimes \phi$\,.
Then $E$ is a normal conditional expectation of $\N$ onto
$\N_0\bar{\otimes} 1$ and $ E(\N^\theta)\supseteq E(R_\infty
\bar{\otimes} 1)=R_\infty \bar{\otimes} 1$\,. Since $\omega$ is
$\theta^{(1)}$ invariant, it follows that $E(\N^\theta)\subseteq
\N_0^{\theta^{(0)}}=R_\infty \bar{\otimes} 1$\,, and therefore
$E_0:=E_0\vert_{\N^\theta}$ is a normal conditional expectation of
$\N^\theta$ onto $R_\infty \bar{\otimes} 1$\,. Hence $\N^\theta$
satisfies condition $(1)$ in Theorem \ref{th8}.

We next prove that $\M\simeq \N^\theta$\,, which will complete the
proof of the implication $(3)\Rightarrow (1)$ in Theorem \ref{th8}
in the hyperfinite case. Note that by (\ref{eq22332222243}), $(\N,
\tau, \theta)$ satisfies condition $(i)$ in Theorem XII.1.1. of
\cite{Tk}, Vol. II. Therefore
\[ Z(\M_1)={Z(\N)}^\theta\,, \]
where $\M_1:=\N\rtimes_{\theta} \mathbb{R}$\,. But since $(Z(\N),
\theta\vert_{Z(\N)})\simeq (A, \theta^{(1)})$\,, the latter being
ergodic, we actually have $Z(\M_1)=\mathbb{C} 1$\,, i.e., $\M_1$ is
a factor. Moreover, by Proposition X.2.6 and Lemma XII.1.2. in
\cite{Tk}, Vol. II, \[ (\N, \mathbb{R}, \theta)\simeq
(\M_2\rtimes_\alpha \mathbb{R}, \mathbb{R}, \widetilde{\alpha})\,,
\] for some covariant system $(\M_2\,, \mathbb{R}\,, \alpha)$\,. By
Theorem X.2.3. in \cite{Tk}, Vol. II, it follows that $\M_2\simeq
\N^\theta$ and
\begin{equation}\label{eq77887766554423232}
\M_1=\N\rtimes_\theta \mathbb{R}\simeq \M_2\bar{\otimes} {\mathcal
B}(L^2(\mathbb{R}))\,.
\end{equation}
Since $\M_2$ is also a factor, and $R_\infty\bar{\otimes} 1$ is the
range of a normal conditional expectation $E_0:\N^\theta\rightarrow
R_\infty\bar{\otimes} 1$\,, it follows by Lemma \ref{lem767654} that
$\M_2\simeq \N^\theta$ is of type III. Hence, by
(\ref{eq77887766554423232}), \[ \N^\theta\simeq \M_2\simeq \M_1\,.
\] By Theorem XII.1.1. in \cite{Tk}, Vol. II, $(\N, \mathbb{R},
\theta)$ is isomorphic to the "noncommutative flow of weights" of
$\M_1$ (in the sense of Def. XII.1.3 in \cite{Tk}, Vol. II) and
hence $(Z(\N), \theta\vert_{Z(\N)})\simeq (A, \theta^{(1)})$ is
isomorphic to the flow of weights for $\M_1$\,. Hence $\M$ and
$\M_1$ are hyperfinite type III-factors with isomorphic flow of
weights, so by {\bf 1.7}, {\bf 1.8} and {\bf 1.9} in the
Introduction, $\M\simeq \M_1$\,, and thus $\M\simeq \N^\theta$\,.
The proof is complete.
\end{rem}

\begin{cor}\label{cor9909}
Let $\M$ be a type III-factor with separable predual. If there
exists a type III$_1$-factor which embeds into $\M$ as the range of
a normal faithful conditional expectation, then there exists a von
Neumann algebra embedding $j: R_\infty\hookrightarrow\M$\,, with a
normal faithful conditional expectation $\widetilde{E}:\M\rightarrow
j(R_\infty)$\,.
\end{cor}

\begin{proof}
Suppose there is a von Neumann algebra embedding $i:{\mathcal
R}\hookrightarrow \M$\,, where ${\mathcal R}$ is a type III$_1$-
factor, with a normal faithful conditional expectation
$E:\M\rightarrow i({\mathcal R})$\,. By the proof of the implication
$1)\Rightarrow 3)$ in Theorem \ref{th8}, we deduce the existence of
an invariant normal state on the flow of weights for $\M$\,. An
application of Theorem \ref{th8} yields the assertion.
\end{proof}
The following lemma is well-known. For completeness, we include a
proof.

\begin{lemma}\label{lem76767632}
Let $\M$ be a hyperfinite factor with separable predual. Then
\[ M\bar{\otimes} R_\infty\cong R_\infty\,. \]
\end{lemma}

\begin{proof}
Let $0< \lambda< 1$\,. Choose $0< \mu< \infty$ such that
$\frac{\log{\lambda}}{\log{\mu}}\notin \mathbb{Q}$\,. Then
$R_\infty\cong R_\mu \bar{\otimes} R_\lambda$ and therefore
$R_\infty\bar{\otimes} R_\lambda\cong R_\infty$\,. We deduce that
$(M\bar{\otimes} R_\infty)\bar{\otimes} R_\lambda\cong
M\bar{\otimes} R_\infty$\,. By the definition of the Araki and Woods
$r_\infty$-invariant of a factor (see Section 3.6 in \cite{Co4}), it
follows that $\lambda \in r_\infty(M\bar{\otimes} R_\infty)$\,. By
\cite{Co4}, Theorem 3.6.1, we conclude that $\lambda\in
S(M\bar{\otimes} R_\infty)$\,, where $S$ is Connes' $S$-invariant of
a factor. Hence $M\bar{\otimes} R_\infty$ is a type III$_1$-factor.
It is also hyperfinite, and thus the assertion follows (see {\bf
1.9}. in the Introduction).
\end{proof}

\noindent{\em Proof of Theorem \ref{th2}}. This follows now
immediately from the equivalence $1)\Leftrightarrow 4)$ in Theorem
\ref{th7}, together with Theorem \ref{th8} and Lemma
\ref{lem555543212} above, and the fact that if $\M$ is an injective
factor with separable predual, then the von Neumann algebra tensor
product $\M\bar{\otimes} R_\infty$ is (isomorphic to) $R_\infty$
(cf. Lemma \ref{lem76767632}), which implies that there exists a von
Neumann algebra embedding $i:\M\hookrightarrow R_\infty$ with a
normal conditional expectation $E:R_\infty\rightarrow i(\M)$\,.\qed

If $\M$ is any factor (not necessarily hyperfinite) of type
III$_\lambda$\,, where $0< \lambda\leq 1$\,, then there always
exists a normal invariant state on the flow of weights for $\M$\,.
Using results of Haagerup and Winsl{\o}w (cf. \cite{HW}; see also
\cite{HRS}, Theorem 6.2), we exhibit in the following an uncountable
family of mutually non-isomorphic (in the von Neumann algebras
sense) hyperfinite type III$_0$-factors which admit a normal
(faithful) invariant state on their flow of weights:
\begin{exam}\label{exp768}\rm
Let $G$ be a dense, countable subgroup of $\mathbb{R}$\,. Further,
let $\phi$ be a normal, faithful state on $R_\infty$ and set
$\N_G:=R_\infty\rtimes_{\alpha}G$\,, where $\alpha:G\rightarrow
\text{Aut}(\M)$ is the restriction of the modular automorphism group
$(\sigma_t^{\phi})_{t\in \mathbb{R}}$ to $G$. Then $\N_G$ is an
injective type III$_0$-factor. Moreover, $T(\N_G)=G$\,, where $T$ is
Connes $T$-invariant. In particular, if $G\neq G'$\,, then $\N_G$
and $\N_{G'}$ are not von Neumann algebras isomorphic. It is easily
checked that there are uncountably many dense, countable subgroups
of $\mathbb{R}$\,. Since $\N_G$ is a crossed product of $R_\infty$
by a discrete group, there exists an embedding $i$ of $R_\infty$
into $\N_G$ with a normal, faithful conditional expectation
$E:\N_G\rightarrow i(R_\infty)$\,. By Theorem \ref{th8}, we deduce
the existence of a normal invariant state on the flow of weights for
$\N_G$\,. Note that $(\N_G)_*$ is cb-isomorphic to $(R_\infty)_*$\,,
as shown in the proof of Theorem 6.2 in \cite{HRS}.
\end{exam}

We end this section with the following results concerning the
non-hyperfinite case.

\begin{prop}\label{pr10}
Let $\M$ be any type III-factor with separable predual $\M_*$\,,
such that there exists an invariant normal state on the flow of
weights for $\M$\,. If $\N$ is any semifinite von Neumann algebra
with separable predual $\N_*$\,, then $\M_*$ is not cb-isomorphic to
a subspace of $\N_*$\,. In particular, $\M_*$ and $\N_*$ are not
cb-isomorphic.
\end{prop}

\begin{proof} By Theorem \ref{th8} it follows that there exists an embedding
$i:R_\infty\rightarrow \M$ with a normal faithful conditional
expectation $E:\M\rightarrow i(R_\infty)$\,. This yields a
cb-embedding $E_*:(R_\infty)_*\hookrightarrow \M_*$ such that
$E_*((R_\infty)_*)$ is cb-complemented into $\M_*$\,.

Suppose by contradiction that $\M_*$ is cb-isomorphic to a subspace
of $\N_*$\,. This gives rise to a a cb-embedding of $(R_\infty)_*$
into $\M_*$\,. The crucial point is now the fact that the operator
Hilbert space $OH$ is cb-isomorphic to a subspace of
$(R_\infty)_*$\,, as proved by Junge \cite{Ju}. A different proof
that yields to the improved cb-isomorphism constant $\leq \sqrt{2}$
has been recently obtained by the authors (see \cite{HM}). However,
since $\N$ is semifinite, $OH$ does not cb-embed into the predual
$\N_*$\,, as shown by Pisier \cite{Pi3}\,. This leads to a
contradiction, and the proof is complete.
\end{proof}
Examples of non-injective type III$_0$-factors with an invariant
normal state on their flow of weights can be obtained as follows:

\begin{exam}
\label{ex6}\rm Given a dense countable subgroup $G$ of $\mathbb{R}$,
let $\N_G$ be the corresponding injective type III$_0$-factor
constructed in Example \ref{exp768}\,. Let $\N$ be a semifinite
non-injective factor, and set $\M_G:=\N_G\bar{\otimes} \N$\,. Then
$\M_G$ is a type III$_0$-factor, since by Corollary 3.2.8 in
\cite{Co4}, $S(\N_G\bar{\otimes} \N)=S(\N_G)$\,. Furthermore, $\M_G$
is clearly non-injective, and, moreover, by Theorem \ref{th8} it
admits an invariant, normal state on the flow of weights.
\end{exam}

\begin{rem}\label{rem2223432}\rm
The existence of an invariant normal state on the flow of weights
appears in a different context in Connes's paper \cite{Co5}. Here a
certain class class of foliated 3-manifolds $(V, F)$ is considered,
and it is proved that if the Godbillon-Vey invariant of such a
foliated manifold is non-zero, then the von Neumann algebra $\M$
associated to $(V, F)$ must have an invariant normal state on its
flow of weights (see \cite{Co5}, Theorem 0.3. and Theorem 7.14.)
\end{rem}

\section{CB-isomorphism classes of preduals of injective type III$_0$-factors}
\setcounter{equation}{0}

In this section we will prove the following:

\begin{theorem}\label{th67678}
For every $0\leq t< 2$\,, there exists a non-transitive ergodic flow
$(A^{(t)}, (\theta_s^{(t)})_{s\in \mathbb{R}})$ with separable
predual $A_*^{(t)}$\,, such that for all normal states $\omega$ on
$A^{(t)}$ we have
\[ \lim\limits_{n\rightarrow \infty}\|\omega\circ
\theta_{2^n}^{(t)}-\omega\|=t\,. \]
\end{theorem}

As a consequence we obtain:

\begin{theorem}\label{th56747}
If $({\M}^{(t)})_{0\leq t< 2}$ are the hyperfinite type III$_0$-
factors with $(A^{(t)}, (\theta_s^{(t)})_{s\in \mathbb{R}})$ from
Theorem \ref{th67678} as flow of weights, then for any $0\leq t_1<
t_2< 2$\,, the predual ${\M}^{(t_1)}_*$ of ${\M}^{(t_1)}$ is not
cb-isomorphic to a cb-complemented subspace of the predual
${\M}^{(t_2)}_*$ of ${\M}^{(t_2)}$\,.

In particular, $({\M}^{(t)}_*)_{0\leq t< 2}$ is a family of mutually
not cb-isomorphic preduals of hyperfinite type III$_0$ factors.
\end{theorem}

\begin{proof}
Let $0\leq t_1< t_2< 2$\,. Assume by contradiction that
${\M}^{(t_1)}_*$ is cb-isomorphic to a cb-complemented subspace of
${\M}^{(t_2)}_*$\,. By Proposition \ref{thprop4} we can then embed
${\M}^{(t_1)}$ into ${\M}^{(t_2)}$ as a von Neumann algebra with a
normal faithful conditional expectation $E: {\M}^{(t_2)}\rightarrow
{\M}^{(t_1)}$ onto. Let $\phi_1$ be a normal faithful state on
${\M}^{(t_1)}$ and set $\phi_2:=\phi_1\circ E$\,. Then $\phi_2$ is a
normal faithful state on ${\M}^{(t_2)}$ satisfying (\ref{eq4567}).
As explained in the proof of Theorem \ref{th8}, this ensures the
existence of an embedding of
${\N}^{(t_1)}:={\M}^{(t_1)}\rtimes_{\sigma^{\phi_1}} \mathbb{R}$
into ${\N}^{(t_2)}:={\M}^{(t_2)}\rtimes_{\sigma^{\phi_2}}
\mathbb{R}$\,, with a normal faithful conditional expectation
$\widetilde{E}: {\N}^{(t_2)}\rightarrow {\N}^{(t_1)}$\,, satisfying
\begin{equation}\label{eq77777779774}
\widetilde{E}\vert_{\pi_2(\M^{(t_2)})}=E \end{equation} and
(\ref{eq445532211})\,, where $\pi_i(\M^{(t_i)})$ and
$(\lambda_i(t))_{t\in \mathbb{R}}$ denote the generators of
$\N^{(t_i)}$\,, $i=1, 2$. By (\ref{eq333334333339}) and
(\ref{eq3333343333392}), we then infer that
\begin{equation*}
\widetilde{E} \circ \theta_s^{(t_2)}=\theta_s^{(t_1)}\circ
\widetilde{E}\,, \quad s\in \mathbb{R}\,.
\end{equation*}
Furthermore, by the bimodule property of conditional expectations,
we have as in the proof of (\ref{eq4568}) that
\[ \widetilde{E}(A^{(t_2)})\subseteq A^{(t_1)}\,. \]
Set $S:=\widetilde{E}\vert_{A^{(t_2)}}$\,. Then $S:
A^{(t_2)}\rightarrow A^{(t_1)}$ is a unital, normal, positive
mapping satisfying
\begin{equation}\label{eq343422327}
S\circ \theta_s^{(t_2)}=\theta_s^{(t_1)}\circ S\,, \quad  s\in
\mathbb{R}\,.
\end{equation}
Choose a normal state $\omega_1$ on $A^{(t_1)}$ and set
$\omega_2:=\omega_1\circ S$\,. Then $\omega_2\in
S_{\text{nor}}(A^{(t_2)})$ and by (\ref{eq343422327}) we get
\[ (\omega_2\circ \theta_s^{(t_2)}-\omega_2)=(\omega_1\circ
\theta_s^{(t_1)}-\omega_1)\circ S\,. \] By positivity we infer that
$\|S\|\leq 1$ and therefore
\[ \|\omega_2\circ \theta_s^{(t_2)}-\omega_2\|\leq \|\omega_1\circ
\theta_s^{(t_1)}-\omega_1\|\,. \] Let $s:=2^n$\,, $n\geq 1$\,, and
pass to the limit as $n\rightarrow \infty$\,. We infer that $t_2\leq
t_1$\,, which is a contradiction.
\end{proof}

The proof of Theorem \ref {th67678} will be achieved in several
steps. Recall first the action of $\mathbb{Z}$ on $\Omega=\{0,
1\}^\infty$ by the {\em dyadic odometer} transformation (cf.
\cite{Tk}, Vol. III, Definition 3.24). Define $g: \Omega\rightarrow
\Omega$ by:
\begin{eqnarray*}
g(0, x_2, x_3, x_4, \ldots) &=& (1, x_2, x_3, x_4, \ldots)\\
g(1, 0, x_3, x_4, \ldots) &=& (0, 1, x_3, x_4, \ldots)\\
\vdots &&\\
g(1, 1, \ldots, 1, 0, x_{n+1}, \ldots) &=& (0, 0, \ldots, 0, 1,
x_{n+1}, \ldots)\\
\vdots &&\\
g(1, 1, 1, \ldots) &=& (0, 0, 0, \ldots)
\end{eqnarray*}
Take $0< a\leq \frac12$\,. Define a measure $\nu_a$ on $\Omega$ by
\[ \nu_a:=\otimes_{n=1}^\infty \mu_a\,, \quad \text{where}
\,\,\,\,\mu_a:=a\delta_0 + (1-a)\delta_1\,. \] It is easy to check
that $g$ preserves the measure class of $\nu_a$\,, i.e., the image
measure $g(\nu_a)$ has the same null-sets as $\nu_a$\,. Therefore
$g$ induces an automorphism $\sigma$ of $L^\infty(\Omega, \nu_a)$ by
\begin{equation}\label{eq88998877667788}\sigma(f)(x)=f(g^{-1}x)\,, \quad f\in L^\infty(\Omega,
\nu_a)\,, x\in \Omega\,. \end{equation} It is well-known that
$\sigma$ is ergodic, i.e., if $f\in L^\infty(\Omega, \nu_a)$
satisfies $\sigma(f)=f$ ($\nu_a$-a.e.), then $f$ is $\nu_a$-a.e.
equal to a constant function. This can be seen by observing that
$(g^n)_{n\in \mathbb{Z}}$ have the same orbits in $\Omega$ (up to
null-sets) as the natural action of
${\mathbb{Z}}_2^{(\infty)}=\{(x_1, x_2, x_3 \ldots )\in \Omega;
x_i\neq 0 \,\,\text{eventually, as} \,\, i\rightarrow \infty\}$ on
$\Omega$\,, generated by $(\rho_k)_{k=1}^\infty$\,, where $\rho_k$
changes $x_k$ to $1-x_k$ in $x=(x_1\,, x_2\,, \ldots)\in \Omega$ and
leaves the remaining coordinates of $x$ unchanged. The induced
action on $L^\infty(\Omega, \nu_a)$
\[ \alpha_\rho(f)(x)=f(\rho^{-1}x)\,, \quad f\in L^\infty(\Omega,
\nu_a)\,, \rho\in {\mathbb{Z}}_2^{(\infty)}\,, x\in \Omega \] is
ergodic, because the crossed-product $L^\infty(\Omega,
\nu_a)\rtimes_\alpha {\mathbb{Z}}_2^{(\infty)}$ is a factor (cf.
\cite{Kr2}, Introduction). Therefore, $\sigma$ is ergodic, as well.

In the following, we will use the symbol $\nu_a$ also to denote the
normal state on $L^\infty(\Omega, \nu_a)$ given by \[
\nu_a(f)=\int_\Omega f \,d\nu_a\,. \] With this notation, we have
\begin{lemma}\label{lem1112232} The following equality holds:
\begin{equation}\label{eq333335}
\|\nu_a\circ \sigma-\nu_a\|=2-4a\,.
\end{equation}
\end{lemma}

\begin{proof}
Note that
\[ (\nu_a\circ \sigma)(f)=\int_\Omega f(g^{-1}x)
\,d\nu_a(x)=\int_\Omega f \,dg^{-1}(\nu_a)\,, \quad f\in
L^\infty(\Omega, \nu_a)\,, \] where $g^{-1}(\nu_a)$ is the image
measure of $\nu_a$ by the map $g^{-1}$\,. Hence
\[ \|\nu_a\circ \sigma-\nu_a\|_{L^\infty(\Omega,
\nu_a)_*}=\left\|\frac{dg^{-1}(\nu_a)}{d\nu_a}-1\right\|_{L^1(\Omega,
\nu_a)}\,. \] For every $n\geq 1$\,, let $K_n$ be the set of
elements in $\Omega$ of the form \[ (\underbrace{1, 1, \ldots,
1}_{n-1\,\, \text{times}}, 0\,, x_{n+1}\,, x_{n+2}\,, \ldots)\,, \]
where $x_j\in \{0, 1\}$ for $j\geq n+1$\,, and put $k_0=(1, 1, 1,
\ldots )$\,. Then $\Omega=\left(\cup_{n=1}^\infty K_n\right)\cup
\{k_0\}$ (disjoint union). By the definition of $g$ and $\nu_a$ it
is clear that for every Borel set $E$ in $K_n$,
\[
\nu_a(gE)=\frac{1-a}{a}\left(\frac{a}{1-a}\right)^{n-1}\nu_a(E)\,.
\]
Since $g^{-1}(\nu_a)(E)=\nu_a(gE)$\,, it follows that
\begin{equation}\label{eq333444433344}
\frac{dg^{-1}(\nu_a)}{d\nu_a}(x)=\left(\frac{a}{1-a}\right)^{n-2}\,,
\quad x\in K_n\,. \end{equation}  Since $1-a\geq a$\,, it follows
that
\begin{eqnarray*} \left\|\frac{dg^{-1}(\nu_a)}{d\nu_a}-1\right\|_{L^1(\Omega,
\nu_a)}&=&\sum\limits_{n=1}^\infty \left\vert
\left(\frac{a}{1-a}\right)^{n-2}-1\right\vert \nu_a(K_n)\\&=&
\left(\frac{1-a}{a}-1\right)a+\sum\limits_{n=3}^\infty
\left(1-\left(\frac{a}{1-a}\right)^{n-2}\right)a(1-a)^{n-1}\\&=&
2-4a\,.
\end{eqnarray*} The last equality can, of course, be obtained by
summation of the infinite sum, but it can be obtained more easily by
observing that, since \[ \int_\Omega \frac{dg^{-1}(\nu_a)}{d\nu_a}
\,d\nu_a=1\,, \] the positive and negative parts of
$\frac{dg^{-1}(\nu_a)}{d\nu_a}-1$ have the same $L^1$-norm, and
therefore
\[ \sum\limits_{n=3}^\infty
\left(1-\left(\frac{a}{1-a}\right)^{n-2}\right)a(1-a)^{n-1}=\left(\frac{1-a}{a}-1\right)a=1-2a\,.
\]
The proof is complete.
\end{proof}

Note that the action of $g$ on $\Omega$ can be considered as binary
addition of $(1, 0, 0, \ldots)$ and $(x_1\,, x_2\,, \ldots)\in
\Omega$ with carry over to the right. More generally, if $k\in
\mathbb{N}$ has the binary representation \[ k=k_1+{k_2}2+
{k_3}2^2+\ldots + {k_m}2^{m-1}\,, \] then the action of $g^k$ on
$\Omega$ is given by binary addition of $(k_1\,, k_2\,, \ldots \,,
k_m\,, 0, 0, \ldots)$ and $(x_1\,, x_2\,, \ldots)$ with carry over
to the right. In particular, if $k=2^n$\,, for some $n\in
\mathbb{N}$\,, then $m=n+1$ and
\[ (k_1\,, k_2\,, \ldots \,,
k_m\,, 0, 0, \ldots)=(\underbrace{0, 0, \ldots, 0}_{n\,\,
\text{times}}, 1, 0, \ldots)\,. \] Hence $g^{2^n}((x_1, \ldots, x_n,
x_{n+1}, x_{n+2}\, \ldots))=(x_1, \ldots , x_n, g((x_{n+1}, x_{n+2},
\ldots)))$\,. This also implies that for all positive integers $n$
we have $\sigma^{2^n}=(\otimes_{i=1}^n {\text{Id}})\otimes
\sigma$\,.

\begin{prop}\label{prop98761}
For all $\phi\in S_{\text{nor}}(L^{\infty}(\Omega, \nu_a))$,
\begin{equation}\label{eq3222212222}
\lim\limits_{n\rightarrow \infty} \|\phi\circ
\sigma^{2^n}-\phi\|=2-4a\,.
\end{equation}
\end{prop}

\begin{proof}
Let $\phi\in S_{\text{nor}}(L^{\infty}(\Omega, \nu_a))$\,. For all
$n\geq 1$\,, let $A_n:=\otimes_{k=1}^n l^\infty \{0, 1\}\otimes 1$
and set $\omega_n:=\phi\vert_{A_n}$\,. Then
\[ \phi\circ E_n=\omega_n\otimes \left(\otimes
_{k={n+1}}^\infty \mu_a\right)\,, \] where $E_n$ denotes the natural
conditional expectation of $L^{\infty}(\Omega, \nu_a)$ onto $A_n$\,.
By standard infinite tensor product theory, it follows that for all
$\psi\in S_{\text{nor}}(L^\infty(\Omega, \nu_a))$\,, we have
\[ \lim\limits_{n\rightarrow \infty} \|\psi\circ E_n-\psi\|=0\,. \]
Hence it suffices to prove (\ref{eq3222212222}) for all $n\geq 1$
and all $\phi\in S_{\text{nor}}(L^{\infty}(\Omega, \nu_a))$ of the
form
\[ \phi:=\omega\otimes \left(\otimes
_{k={n+1}}^\infty \mu_a\right)\,, \quad \omega\in
S_{\text{nor}}(A_n)\,. \] Fix $n\geq 1$ and consider $\phi\in
S_{\text{nor}}(L^{\infty}(\Omega, \nu_a))$ of this form. Then, for
all $m\in \mathbb{N}$ with $m> n$, we have
\[ \phi=\omega\otimes \underbrace{\mu_a\otimes \ldots \otimes
\mu_a}_{m-n\,\, \text{times}}\otimes \left(\otimes_{k={n+1}}^\infty
\mu_a\right)\,. \] Since $\left(\otimes_{k={n+1}}^\infty
\mu_a\right)=\nu_a$\,, we conclude by previous considerations that
\[ \phi\circ \sigma^{2^m}=\omega \otimes \underbrace{\mu_a\otimes
\ldots \otimes \mu_a}_{m-n\,\, \text{times}}\otimes (\nu_a\circ
\sigma)\,. \] By Lemma \ref{lem1112232} we deduce that
\begin{eqnarray*}
\|\phi\circ \sigma^{2^m}-\phi\|&=&\|\omega \otimes
\underbrace{\mu_a\otimes \ldots \otimes \mu_a}_{m-n\,\,
\text{times}}\otimes (\nu_a\circ \sigma-\nu_a)\|\\
&=&\|\nu_a\circ \sigma-\nu_a\|\\&=&2-4a\,.
\end{eqnarray*}
Hence $\lim\limits_{m\rightarrow \infty} \|\phi\circ
\sigma^{2^m}-\phi\|=2-4a$\,, and the proof is complete.
\end{proof}

{\em Proof of Theorem \ref{th67678}}. Let $(g, \Omega, \nu_a)$ be
the transformation space as above. In the following we will
construct the associated flow under the constant ceiling function
$1$. This is a special case of Krieger's construction of the flow
under a ceiling function $\phi$ given in \cite{Kr1}, p. 46-47.

Set $\widetilde{\Omega}:=\Omega\times [0, 1)$ and
$\widetilde{\nu_a}=\nu_a\times dx$\,, where $dx$ is the Lebesgue
measure on $[0, 1)$\,. Given $s\in \mathbb{R}$\,, define
\begin{equation}\label{eq3332222279}
\tilde{g}_s(x, y):=(g^n(x), t')\,, \quad x\in \Omega\,, 0\leq y<
1\,,
\end{equation}
where $s+y=n+y'$\,, with $n\in \mathbb{Z}$ and $0\leq y'< 1$\,.

Then $((\tilde{g}_s)_{s\in \mathbb{R}}, \widetilde{\Omega},
\widetilde{\nu_a})$ is a one-parameter group of Borel-measurable
actions on $\widetilde{\Omega}$ which preserve the measure class of
$\widetilde{\nu_a}$\,. Define now $\widetilde{\sigma}$ as the
corresponding action on $L^{\infty}(\widetilde{\Omega},
\widetilde{\nu_a})$\,, i.e., for all $s\in \mathbb{R}$ and all $f\in
L^{\infty}(\widetilde{\Omega}, \widetilde{\nu_a})$ let
\[ ({\widetilde{\sigma}}_s(f))(z):=f(\widetilde{g}_s^{-1}(z))\,,
\quad z=(x, y)\in \widetilde{\Omega}\,. \] For simplicity of
notation, set $\widetilde{\sigma}:=(\widetilde{\sigma}_s)_{s\in
\mathbb{R}}$\,.

By the remark following Definition 3.1 in \cite{Tk} (Vol. II, p.
385), $\widetilde{\sigma}$ is ergodic. Also, $g$ is non-transitive,
because every orbit $\{g^nx; n\in \mathbb{Z}\}$ is countable, so all
$g$-orbits have $0$ measure. Hence
$\widetilde{g}:=(\widetilde{g}_s)_{s\in \mathbb{R}}$ is
non-transitive, because $\widetilde{g}$-orbits are of the form
$L=L_0\times [0, 1)$\,, where $L\subseteq \Omega$ is an orbit for
$g$, so
\[ \widetilde{\nu_a}(L)=\nu_a(L_0)=0\,. \]
We conclude that $(L^{\infty}(\widetilde{\Omega},
\widetilde{\nu_a}), \widetilde{\sigma})$ is an ergodic and
non-transitive flow. Hence this flow is the "smooth flow of weights"
of a unique (up to von Neumann algebras isomorphism) hyperfinite
factor of type III$_0$ (see {\bf 1.4} in the Introduction).

We claim that for all $\phi\in
S_{\text{nor}}(L^{\infty}(\widetilde{\Omega}, \widetilde{\nu_a}))$,
\begin{equation}\label{eq322221222278}
\lim\limits_{n\rightarrow \infty} \|\phi\circ
{\widetilde{\sigma}}_{2^n}-\phi\|=2-4a\,.
\end{equation}
Indeed, let $k\in \mathbb{Z}$\,. By (\ref{eq3332222279}) we have
\[ \widetilde{g}_k(x, y)=(g^k(x), y)\,, \quad x\in \Omega\,,
0\leq y< 1\,. \] Hence
\begin{equation}\label{eq76767655554321}{\widetilde{\sigma}}_k=\sigma^k\otimes \text{Id}_{L^\infty([0, 1)\,,
dx)}\,.\end{equation} In particular, we deduce for any $n\geq 1$
that ${\widetilde{\sigma}}_{2^n}=\sigma^{2^n}\otimes
\text{Id}_{L^\infty([0, 1)\,, dx)}$\,.

Given a positive integer $m$, denote by $B_m$ the set of functions
which are constant on $\left[\frac{i}{2^m}\,,
\frac{i+1}{2^m}\right)$\,, for all $0\leq i\leq 2^m$\,. Note that
\[
B_m=\text{Span}\left\{q_i:=1_{\left[\frac{i}{2^m}\,,
\frac{i+1}{2^m}\right)}; 0\leq i< 2^m\right\}\,. \] Further, let
$F_m: L^\infty([0, 1))\rightarrow B_m$ be the natural conditional
expectation onto $B_m$ preserving Lebesgue measure $dx$\,. Then the
mapping $\widetilde{F}_m:=\text{Id}_{L^\infty(\Omega, \nu_a)}\otimes
F_m$ is a conditional expectation of $L^\infty(\widetilde{\Omega},
\widetilde{\nu_a})$ onto ${L^\infty(\Omega, \nu_a)}\otimes B_m$
preserving $\widetilde{\nu}_a$\,. Clearly, for all normal states
$\phi$ on $L^{\infty}(\widetilde{\Omega}, \widetilde{\nu}_a)$\,,
\[ \lim\limits_{n\rightarrow \infty} \|\phi\circ
\widetilde{F}_m-\phi\|=0\,. \] Therefore, in order to prove
(\ref{eq322221222278}) it suffices to consider states $\phi$ of the
form \begin{equation}\label{eq767676555543216} \phi:=\omega\circ
\widetilde{F}_m\,, \end{equation} where $m\in \mathbb{N}$
(arbitrarily chosen) and $\omega$ is a normal state on
${L^\infty(\Omega, \nu_a)}\otimes B_m$\,.

Fix now $m\in \mathbb{N}$\,, and let $\phi$ be of the form
(\ref{eq767676555543216}). For $0\leq i< 2^m$, set
\[ \omega_i(f):=\omega(f\otimes q_i)\,, \quad f\in
L^\infty(\Omega, \nu_a)\,. \] Then $\omega_i$ are positive linear
functionals on $L^\infty(\Omega, \nu_a)$\,,
$\sum\limits_{i=1}^{2^m-1} \omega_i(1)=1$, and the $\omega_i$'s
determine $\omega$ uniquely. For any $k\in \mathbb{Z}$ we obtain by
(\ref{eq76767655554321}) that
\[ \|\phi\circ
{\widetilde{\sigma}}_k-\phi\|=\sum\limits_{i=0}^{2^m-1}\|\omega_i\circ
\sigma^k-\omega_i\|\,. \] By Proposition \ref{prop98761} we deduce
that $\lim\limits_{n\rightarrow \infty} \|\chi\circ
\sigma^{2^n}-\chi\|=(2-4a)\|\chi\|$\,, for every positive normal
functional $\chi$ on $L^\infty(\Omega, \nu_a)$\,. Hence
\begin{equation*}
\lim\limits_{n\rightarrow \infty}
\left(\sum\limits_{i=0}^{2^m-1}\|\omega_i\circ
\sigma^{2^n}-\omega_i\|\right)=(2-4a)\sum\limits_{i=0}^{2^m-1}
\|\omega_i\|=(2-4a)\sum\limits_{i=0}^{2^m-1} \omega_i(1)= (2-4a)\,.
\end{equation*}
The proof is complete.\qed

In the following we will compute Connes' $T$-invariant for the
hyperfinite type III$_0$-factors $\M^{(t)}$\,, $0\leq t< 2$
constructed above. Recall that, if $\M$ is a von Neumann algebra
with a normal, faithful state $\phi$\,, then Connes' $T$-invariant
$T(\M)$ defined by
\[ T(\M):=\{\tau\in \mathbb{R}; \, \sigma_\tau^{\phi}\in
\text{Int}(\M)\} \] is independent of $\phi$\,, since for any
normal, faithful state $\psi$ on $\M$\,,
\[ \sigma_\tau^{\psi}(x)=(D\psi: D\phi)_\tau\sigma_\tau^{\phi}(x) (D\psi:
D\phi)_\tau^*\,, \quad \tau\in \mathbb{R}\,, x\in \M\,. \] By
\cite{Tk}, Vol. II, Chap. XII, if $\M$ has flow of weights $(Z(\N),
\theta)$\,, then
\[ T(\M)=\{\tau \in \mathbb{R}; \exists \,u\in Z(\N), u
\,\text{unitary such that}\, \theta_s(u)=e^{i\tau s}u\,, s\in
\mathbb{R}\}\,. \]

\begin{theorem}\label{th9322}
For all $0\leq t< 2$\,,
\begin{equation}\label{eq111122332211}
T(\M^{(t)})=\left\{\frac{2\pi k}{2^n}; k\in \mathbb{Z}\,, n\in
\mathbb{N}\right\}\,.
\end{equation}
\end{theorem}

The proof is based on the following intermediate results:
\begin{lemma}\label{lem3322332211111}
For $0\leq t< 2$\,, let $(A^{(t)}, \theta^{(t)})$ be the flow
constructed in the proof of Theorem \ref{th67678}. Then for all
$f\in A^{(t)}$\,,
\[ \lim\limits_{n\rightarrow \infty} \|\theta_{2^n}^{(t)}
(f)-f\|_2=0\,, \] where the 2-norm is taken with respect to the
measure $\widetilde{\nu}_a=\nu_a\otimes {dx}$\,, defined in the
proof of Theorem \ref{th67678} ($t=2-4a$)\,.
\end{lemma}

\begin{proof}
Let $\Omega$\,, $\nu_a$ and $\sigma$ be as defined above (see
(\ref{eq88998877667788})). Then, for all $f\in L^\infty(\Omega,
\nu_a)$\,, \begin{equation*} \|\sigma(f)\|_2^2= \int_\Omega
|h(g^{-1}x)|^2 d\nu_a(x) =\int_\Omega |h(x)|^2dg^{-1}(\nu_a)(x) =
\int_\Omega |h(x)|^2\frac{dg^{-1}(\nu_a)}{d\nu_a}(x) d\nu_a(x)\,.
\end{equation*}
By (\ref{eq333444433344}), we have
\[ \frac{dg^{-1}(\nu_a)}{d\nu_a}(x)
d\nu_a(x)=\sum\limits_{n=1}^\infty
\left(\frac{a}{1-a}\right)^{n-2}1_{\Omega_n}\,. \] Hence
$\left\|\frac{dg^{-1}(\nu_a)}{d\nu_a}\right\|_\infty\leq
\frac{1-a}{a}$\,, and therefore
\begin{equation}\label{eq77765432}
\|\sigma(f)\|_2^2\leq \frac{1-a}{a}\|f\|_2^2\,.
\end{equation}
Since for all positive integers $n$\,,
$\sigma^{2^n}=\text{Id}_{A_n}\otimes \sigma^{(n)}$\,, where
$\sigma^{(n)}$ is equal to $\sigma$ shifted to
$\otimes_{k=n+1}^\infty l^\infty\{0, 1\}$\,, we get from
(\ref{eq77765432}) that
\[ \|\sigma^{2^n}(f)\|_2\leq
\left(\frac{1-a}{a}\right)^{\frac12}\|f\|_2\,, \quad f\in
L^\infty(\Omega, \nu_a)\,. \] Moreover, since
$\widetilde{\sigma}_{2^n}=\sigma^{2^n}\otimes
\text{Id}_{L^\infty([0, 1)\,, dx)}$\,, it also follows that
\[ \|\widetilde{\sigma}_{2^n}(f)\|_2\leq \left(\frac{1-a}{a}\right)^{\frac12}\|f\|_2\,, \quad f\in
L^\infty(\widetilde{\Omega}_a, \widetilde{\nu}_a)\,. \] Since
$\cup_{n=1}^\infty A_n$ is dense in $L^\infty (\Omega_a, \nu_a)$\,,
it follows that the increasing union $\cup_{n=1}^\infty (A_n\otimes
L^([0, 1)))$ is dense in $L^\infty (\widetilde{\Omega}_a\,,
\widetilde{\nu}_a)$\,. Let now $f\in L^\infty
(\widetilde{\Omega}_a\,, \widetilde{\nu}_a)$ and $\varepsilon> 0$\,.
Choose $n\in \mathbb{N}$ and $g\in A_n\otimes L^([0, 1))$ such that
$\|f-g\|_2< \varepsilon$\,. Since
$\sigma^{2^n}\vert_{A_n}=\text{Id}_{A_n}$\,, we have
$\widetilde{\sigma}_{2^n}(g)=(\sigma^{2^n}\otimes
\text{Id}_{L^\infty([0, 1))})(g)=g$\,. Therefore,
\begin{eqnarray*}
\|\widetilde{\sigma}_{2^n}(f)-f\|_2&\leq &
\|\widetilde{\sigma}_{2^n} (f-g)\|_2+\|g-f\|_2\\
&\leq &
\left(\left(\frac{1-a}{a}\right)^{\frac12}+1\right)\|f-g\|_2\\&\leq
& \left(\left(\frac{1-a}{a}\right)^{\frac12}+1\right)\varepsilon\,.
\end{eqnarray*}
This shows that $\lim\limits_{n\rightarrow \infty}
\|\widetilde{\sigma}_{2^n}(f)-f\|_2=0$\,, and the proof is complete.
\end{proof}

\begin{lemma}\label{lem33223322111115}
Let $\tau\in \mathbb{R}$ such that $\lim\limits_{n\rightarrow
\infty} e^{i\tau 2^n} =1$\,. Then $\tau\in \left\{\frac{2\pi
k}{2^n}; k\in \mathbb{Z}\,, n\in \mathbb{N}\right\}$\,.
\end{lemma}

\begin{proof}
Choose $n_0\in \mathbb{N}$ such that
\begin{equation}\label{eq333111133311}
|e^{i\tau 2^n}-1|< 1\,, \quad n\geq n_0\,.
\end{equation}
Assume further that $e^{i\tau 2^{n_0}}\neq 1$\,. Then by
(\ref{eq333111133311}), it follows that $e^{i\tau
2^{n_0}}=e^{iv}$\,, for some $v\in \left(-\frac{\pi}{3},
\frac{\pi}{3}\right)\setminus \{0\}$\,. Hence there exists $k\in
\mathbb{N}$ such that $\frac{\pi}{3}2^{-k}\leq |v|<
\frac{\pi}{3}2^{1-k}$\,. But then $e^{i\tau 2^{n_0+k}}=e^{iv2^k}$
and $\frac{\pi}{3}\leq 2^k v< \frac{2\pi}{3}$\,. Then $|e^{i\tau
2^{n_0+k}}-1|\geq 1$\,, which contradicts (\ref{eq333111133311}).
Hence $e^{i\tau 2^{n_0}}=1$\,. This yields the conclusion.
\end{proof}

{\em Proof of Theorem \ref{th9322}}: Fix $0\leq t< 2$\,. Let $\tau
\in T(\M_t)$\,. Then there exists a unitary $u\in L^(\Omega, \nu_a)$
such that $\theta_s(u)=e^{i\tau s} (u)$\,, for all $s\in
\mathbb{R}$\,. Since
\[ \theta_{2^n}(u)\overset{\text{s.o.}}{\longrightarrow} u\,, \quad
\text {as}\,\, n\rightarrow \infty\,, \] it follows by Lemma
\ref{lem3322332211111} that $\lim\limits_{n\rightarrow \infty}
e^{i\tau 2^n}=1$\,. By Lemma \ref{lem33223322111115}, we conclude
that $\tau\in \left\{\frac{2\pi k}{2^n}; k\in \mathbb{Z}\,, n\in
\mathbb{N}\right\}$\,.

Conversely, let $\tau\in \left\{\frac{2\pi k}{2^n}; k, n\in
\mathbb{N}\right\}$\,. Then there exists $n\in \mathbb{N}$ such that
$e^{i\tau 2^n}=1$\,. Put, as before, $A_n=\otimes_{k=1}^n l^{\infty}
\{0, 1\} \otimes 1$\,. Note that $\text{dim}(A_n)=2^n$\,. For $0\leq
j\leq 2^{n-1}$\,, put
\[ G_j:=\{(k_1^{(j)}\,, k_2^{(j)}\,, \ldots \,, k_n^{(j)}\,,
x_{n+1}\,, x_{n+2}\,, \ldots ); x_j\in \{0, 1\}
\,\,\,\text{for}\,\,\, j\geq n+1\}\,, \] where
$j=k_1^{(j)}+k_2^{(j)} 2+\ldots +k_n^{(j)} 2^{n-1}$ is the unique
binary representation of $j$ ($k_i^{j}\in\{0, 1\}$). Using again the
fact that the action of $g$ on $\Omega$ is given by the binary
addition of $(1, 0, 0, \ldots)$ and $(x_1\,, x_2\,, \ldots)\in
\Omega$ with carry over to the right, it follows that
\[ g(G_j)=G_{j+1}\,, \quad 0\leq j\leq 2^{n-1}\,, \]
where $G_{2^{n}}=G_0$\,. Hence $p_j:=1_{G_j}$ ($0\leq j\leq
2^{n-1}$) are orthogonal projections in $A_n$ with sum equal to 1,
satisfying $\sigma(p_j)=p_{j+1}$ ($0\leq j\leq 2^{n-1}$)\,, where
indices are calculated modulo $2^n$\,. Set now
\[ u_0:=p_0+e^{-i\tau}p_1+\ldots + e^{-i(2^n-1)\tau}p_{2^n-1}\,. \]
Then $u_0\in A_n$ is unitary and satisfies
$\sigma(u_0)=e^{i\tau}u_0$\,. Next set
\[ u(x, y)=u_0(x)e^{-i\tau y}\,, \quad x\in \Omega\,, y\in [0,
1)\,. \] Then $u$ is a unitary in $L^\infty(\widetilde{\Omega}_a\,,
\widetilde{\nu}_a)$\,. We will check that
\begin{equation}\label{eq66666432166}
\widetilde{\theta}_s u=e^{i\tau s} u\,, \quad s\in \mathbb{R}\,,
\end{equation}
which implies that $\tau \in T(\M^{(t)})$\,. Indeed, for any $s\in
\mathbb{R}$\,,
\[ (\widetilde{\theta}_{-s} u)(x, y)=u(\widetilde{g}_s(x, y))=u(g^nx,
y')\,, \] where $s+y=n+y'$ (integer part and fractional part,
respectively, of $s+y$). Hence
\begin{eqnarray*}
(\widetilde{\theta}_{-s} u)(x, y)\,\,= \,\,u_0(g^n x)e^{-i\tau y'}=
(\sigma^{-n}u_0)(x)e^{-i\tau y'} &=& e^{-in\tau}u_0(x)e^{-i\tau
(s+y-n)}\\ &=& e^{-i\tau s}u_u(x) e^{-i\tau y}\\&=& e^{-i\tau s}
u(x, y)\,.
\end{eqnarray*}
Replace now $s$ by $-s$ to obtain (\ref{eq66666432166}). This
completes the proof of Theorem \ref{th9322}.\qed

\thanks{}

\end{document}